\documentclass[11pt]{amsart}

\usepackage{amssymb}

\begin{document}

\title[Reflection ideals and mappings]{Reflection
ideals and mappings between generic submanifolds in
complex space}
\author[M. S. Baouendi, N. Mir and L. P.
Rothschild]{M. S. Baouendi, Nordine Mir, and Linda
Preiss Rothschild}
\address{Department of Mathematics, 0112, University
of California at San Diego, La Jolla, CA 92093-0112, USA}
\email{sbaouendi@ucsd.edu, lrothschild@ucsd.edu \hfil\break}
\address{Universit\'e de Rouen, Laboratoire de
Math\'ematiques Rapha\"el Salem, UMR 6085 CNRS, 76821
Mont-Saint-Aignan Cedex, France}
\email{Nordine.Mir@univ-rouen.fr \hfil\break}
%\abstract \endabstract
%\keywords \endkeywords
%\subjclass{32H02}
\thanks{2000 {{\em Mathematics Subject Classification.}
32H02}}
\begin{thanks}{The first and third authors are
partially supported by National Science
Foundation grant DMS 98-01258. The second author would
like to thank UCSD for its hospitality during the
completion of this work.}\end{thanks}

\font\lin=msbm10 scaled 900 
\font\lind=msbm10 scaled 800

\def\Nn{{\hbox {\lin N}}}
\def\Nnn{{\hbox {\lind N}}} 

\numberwithin{equation}{section}

\def\Label#1{\label{#1}}
\def\1#1{\ov{#1}}
\def\2#1{\widetilde{#1}}
\def\3#1{\mathcal{#1}}
\def\4#1{\widehat{#1}}

\def\s{s}
\def\k{\kappa}
\def\ov{\overline}
\def\span{\text{\rm span}}
\def\ad{\text{\rm ad }}
\def\tr{\text{\rm tr}}
\def\xo {{x_0}}
\def\Rk{\text{\rm Rk\,}}
\def\sg{\sigma}
\def \emxy{E_{(M,M')}(X,Y)}
\def \semxy{\scrE_{(M,M')}(X,Y)}
\def \jkxy {J^k(X,Y)}
\def \gkxy {G^k(X,Y)}
\def \exy {E(X,Y)}
\def \sexy{\scrE(X,Y)}
\def \hn {holomorphically nondegenerate}
\def\hyp{hypersurface}
\def\prt#1{{\partial \over\partial #1}}
\def\det{{\text{\rm det}}}
\def\wob{{w\over B(z)}}
\def\co{\chi_1}
\def\po{p_0}
\def\fb {\bar f}
\def\gb {\bar g}
\def\Fb {\ov F}
\def\Gb {\ov G}
\def\Hb {\ov H}
\def\zb {\bar z}
\def\wb {\bar w}
\def \qb {\bar Q}
\def \t {\tau}
\def\z{\chi}
\def\w{\tau}
\def\Z{\zeta}

\def \T {\theta}
\def \Th {\Theta}
\def \L {\Lambda}
\def\b {\beta}
\def\a {\alpha}
\def\o {\omega}
\def\l {\lambda}

\def \im{\text{\rm Im }}
\def \re{\text{\rm Re }}
\def \Char{\text{\rm Char }}
\def \supp{\text{\rm supp }}
\def \codim{\text{\rm codim }}
\def \Ht{\text{\rm ht }}
\def \Dt{\text{\rm dt }}
\def \hO{\widehat{\mathcal O}}
\def \cl{\text{\rm cl }}
\def \bR{\mathbb R}
\def \bC{\mathbb C}
\def \C{\mathbb C}
\def \N{\mathbb N}
\def \bL{\mathbb L}
\def \bZ{\mathbb Z}
\def \bN{\mathbb N}
\def \scrF{\mathcal F}
\def \scrK{\mathcal K}
\def \mc #1 {\mathcal {#1}}
\def \scrM{\mathcal M}
\def \cR{\mathcal R}
\def \scrJ{\mathcal J}
\def \scrA{\mathcal A}
\def \scrO{\mathcal O}
\def \scrV{\mathcal V}
\def \scrL{\mathcal L}
\def \scrE{\mathcal E}
\def \hol{\text{\rm hol}}
\def \aut{\text{\rm aut}}
\def \Aut{\text{\rm Aut}}
\def \J{\text{\rm Jac}}
\def\jet#1#2{J^{#1}_{#2}}
\def\gp#1{G^{#1}}
\def\gpo{\gp {2k_0}_0}
\def\emmp {\scrF(M,p;M',p')}
\def\rk{\text{\rm rk\,}}
\def\Orb{\text{\rm Orb\,}}
\def\Exp{\text{\rm Exp\,}}
\def\Span{\text{\rm span\,}}
\def\d{\partial}
\def\D{\3J}
\def\pr{{\rm pr}}

\def\dbl{[\hskip -1pt [}
\def\dbr{]\hskip -1pt]}

\def \CZZ {\C \dbl Z,\zeta \dbr}
\def \D{\text{\rm Der}\,}
\def \Rk{\text{\rm Rk}\,}
\def \ima{\text{\rm im}\,}
\def \I {\mathcal I}

\newtheorem{Thm}{Theorem}[section]
\newtheorem{Def}[Thm]{Definition}
\newtheorem{Cor}[Thm]{Corollary}
\newtheorem{Pro}[Thm]{Proposition}
\newtheorem{Lem}[Thm]{Lemma}
\newtheorem{Rem}[Thm]{Remark}

\maketitle

\tableofcontents

\section{Introduction and main results}\Label{int}
In this paper, we study formal mappings  between smooth
generic submanifolds in $\C^N$ and establish results on
finite determination, convergence and local biholomorphic
and algebraic equivalence. Our finite determination
result gives sufficient conditions to guarantee that a
formal map as above is uniquely determined  by its jet at a
point of a preassigned order. For real-analytic generic
submanifolds, we prove convergence of formal mappings
under appropriate assumptions and also give
natural geometric conditions to assure that if two germs of
such submanifolds are formally equivalent, then they are
necessarily biholomorphically equivalent. If the
submanifolds are moreover real-algebraic, we address the
question of deciding when biholomorphic equivalence implies
algebraic equivalence. In particular, we prove that if two
real-algebraic hypersurfaces in $\C^N$ are biholomorphically
equivalent, then they are in fact algebraically equivalent.
All the results are first proved in the more general context
of ``reflection ideals" associated to formal mappings
between formal as well as real-analytic and real-algebraic
manifolds.

We now give precise definitions in order to state some of
our
main results. Let $p\in \C^N$ and $p'\in \C^{N'}$. A formal
map
$H:(\C^N,p)\to (\C^{N'},p')$ is an $N'$-vector
of formal power series in $Z-p$ with $H(p)=p'$. The map
$H(Z)=(H_1(Z),\ldots,H_{N'}(Z))$ is called {\em finite} if
the quotient ring $\C \dbl Z-p
\dbr/(H(Z))$ is finite dimensional as a vector space over
$\C$, where $(H(Z))$ is the ideal generated by the $H_j(Z)$
in $\C \dbl Z-p \dbr$, $j=1,\ldots,N'$. In the case
$N=N'$,
$H$ is called {\it invertible} if its Jacobian determinant
does not vanish at
$p$.

Recall that a smooth submanifold $M\subset \C^N$ is called
{\em generic} if it is locally defined by the vanishing of
 smooth real-valued functions
$r_1(Z,\bar{Z}),\ldots,r_d(Z,\bar{Z})$ with 
linearly independent complex differentials
$\partial r_1(Z,\bar{Z}),\ldots,
\partial r_{d}(Z,\bar{Z})$. A generic submanifold
$M\subset
\C^N$ is said to be of {\em finite type}  at $p\in M$ in the
sense of Kohn \cite{Kohn} and Bloom-Graham \cite{BG} if the
Lie algebra generated by the $(0,1)$ and $(1,0)$ smooth
vector fields tangent to
$M$ spans the complexified tangent space of $M$ at $p$.

A (holomorphic) formal vector field at $p\in
\C^N$ is given by
$$X=\sum_{k=1}^Na_k(Z)\frac{\partial}{\partial Z_k}$$ 
with $a_k(Z)\in \C \dbl Z-p \dbr$, $k=1,\ldots,N$. If $M$
is a generic submanifold of real codimension $d$ as
above, and $r_1,\ldots,r_d$ are smooth real-valued defining
functions of $M$ near $p\in M$, we denote by $\rho
(Z,\bar{Z})=(\rho_1(Z,\bar{Z}),\ldots,\rho_d(Z,\bar{Z}))$
the Taylor series of $r_1,\ldots,r_d$ at $p$ considered as
formal power series in $Z-p$ and $\bar{Z}-\bar{p}$. A
holomorphic formal vector field
$X$ at $p\in M$ is called {\em tangent to} $M$ if 
$$(X\rho)(Z,\bar{Z})=c(Z,\bar{Z})\rho(Z,\bar{Z}),$$
where $c(Z,\bar{Z})$ is a $d\times d$ matrix with entries
in $\C \dbl Z-p,\bar{Z}-\bar{p} \dbr$.
Following Stanton \cite{St}, we say that the submanifold
$M$ is {\em holomorphically nondegenerate}
at
$p\in M$ if there is no nontrivial formal holomorphic
vector field at $p$ tangent to $M$ (see \cite{BER}, \S
11.7).

Let $M\subset \C^N$ and $M'\subset \C^{N'}$ be smooth
generic submanifolds of codimension $d$ and $d'$ through $p$
and
$p'$ respectively and 
$H:(\C^N,p)\to (\C^{N'},p')$ a formal map. We say that
$H$ maps $M$ into $M'$ and write $H(M)\subset M'$ if 
$$\rho'(H(Z),\overline{H(Z)})=a(Z,\bar{Z})\rho(Z,\bar{Z}),$$
where $\rho (Z,\bar{Z})$ is the $d$-vector valued formal
power series defined as above for $(M,p)$,
$\rho'(Z',\bar{Z}')$ is the $d'$-vector valued corresponding
series for
$(M',p')$, and $a(Z,\bar{Z})$ is a $d'\times d$ matrix with
entries in $\C \dbl Z-p,\bar{Z}-\bar{p} \dbr$.

We are now ready to state some of the main results of this
paper. We will discuss
previous related work towards the end of this introduction.
Our first two results deal with finite determination of
formal mappings between smooth generic submanifolds in
$\C^N$, as well as convergence of such mappings when the
submanifolds are real-analytic.

\begin{Thm}\Label{fd}
Let $M,M'\subset \C^N$ be smooth generic submanifolds of
the same dimension through $p$ and $p'$ respectively.
Assume that
$M$ is of finite type at $p$ and that $M'$ is
holomorphically nondegenerate at $p'$. Let $H^0:(\C^N,p)\to
(\C^N,p')$ be a formal finite map sending $M$ into $M'$.
Then there exists an integer $K$ such that if
$H:(\C^N,p)\to (\C^N,p')$ is another formal map sending $M$
into $M'$ with
$$\partial^{\alpha}H(p)=\partial^{\alpha}H^0(p),\quad
|\alpha|\leq K,$$
it follows that $H=H^0$.
\end{Thm}

We should mention that Theorem \ref{fd} is new even
if $H$ is assumed to be holomorphic and $M,M'$ are real-analytic. As
an application of Theorem \ref{fd}, it follows for example that if
$h^0:(M,p)\to (M',p')$ is a germ of a smooth CR diffeomorphism 
with
$M$ and $M'$ satisfying the assumptions of Theorem
\ref{fd} and if $h:(M,p)\to (M',p')$ is another
smooth CR map whose Taylor polynomial of order $K$ at $p$
agrees with that of $h^0$, then necessarily the
entire Taylor series at $p$ of $h$ and $h^0$ are the same.

\begin{Thm}\Label{conv}
Let $M,M'\subset \C^N$ be real-analytic generic
submanifolds of the same dimension through $p$ and $p'$
respectively. Assume that $M$ is of finite type at $p$ and
that $M'$ is holomorphically nondegenerate at $p'$. Then
any formal finite map $H:(\C^N,p)\to (\C^N,p')$ sending $M$
into $M'$ is necessarily convergent.
\end{Thm}
It is worth mentioning that the holomorphic nondegeneracy
condition in Theorems \ref{fd} and \ref{conv} is necessary for the
conclusions of those theorems to hold (see \S \ref{open} for comments
and details).

We say that two germs $(M,p)$ and $(M',p')$ of smooth
generic submanifolds in $\C^N$ of the same dimension are {\em
formally equivalent} if there exists a formal invertible map
$H:(\C^N,p)\to (\C^N,p')$ sending $M$ into $M'$. If $M$ and $M'$ are
real-analytic and the invertible map
$H$ can be chosen to be convergent, we say that $(M,p)$ and
$(M',p')$ are {\em biholomorphically equivalent}. Two
formal mappings 
$H,\check H: (\C^N,p)\to (\C^{N'},p')$ are said to {\em agree up order
$\kappa$}, where 
$\kappa$ is a positive integer,  if their
Taylor series at $p$ agree up to order $\kappa$. The
following theorem may be viewed as an approximation result
for formal mappings between 
real-analytic generic submanifolds by convergent mappings,
in the spirit of Artin's approximation theorem
\cite{artin1}.

\begin{Thm}\Label{approx1}
Let $(M,p)$ and $(M',p')$ be two germs of real-analytic
generic submanifolds in $\C^N$ of the same dimension with
$M$ of finite type at $p$. If $H:(\C^N,p)\to (\C^N,p')$ is a
formal finite map sending $M$ into $M'$ and if $\kappa$
is a positive integer, then there exists a convergent map
$H^{\kappa}:(\C^N,p)\to (\C^N,p')$ which sends
$M$ into $M'$ and agrees with $H$ up to order
$\kappa$.
\end{Thm}

We should point out that the assumptions of Theorem
\ref{approx1} do not imply that the given formal map $H$
is itself convergent. The following, which is an immediate
corollary of Theorem
\ref{approx1}, concerns formal and biholomorphic
equivalence.
\begin{Cor}\Label{fb}
Let $(M,p)$ and $(M',p')$ be two germs of real-analytic
generic submanifolds in
$\C^N$ of the same dimension with $M$ of finite type at $p$.
Then
$(M,p)$ and $(M',p')$ are formally equivalent if and only
if
they are biholomorphically equivalent.
\end{Cor}
A convergent mapping $H:(\C^N,p)\to (\C^{N'},p')$ is called
{\em algebraic} if each of its components satisfies a
non-trivial polynomial equation with holomorphic
polynomial coefficients. A germ of a real-analytic generic
submanifold $(M,p)$ in $\C^N$ is called {\em real-algebraic}
if it is contained in a real-algebraic subset of $\C^N$ of
the same real dimension as that of $M$. We say that two
germs
$(M,p)$ and $(M',p')$ of real-algebraic generic
submanifolds of
$\C^N$ of the same dimension are {\em algebraically
equivalent} if there is a germ of an invertible algebraic
map
$H:(\C^N,p)\to (\C^{N},p')$ sending $M$ into $M'$. The
following theorem can be viewed as an approximation
result for local holomorphic mappings between 
real-algebraic generic submanifolds by algebraic mappings.

\begin{Thm}\Label{approx2}
Let $M,M'\subset \C^N$ be two real-algebraic generic
submanifolds of the same dimension. Assume that $M$ is
connected and of finite type at some point. Let
$p\in M$, $p'\in M'$ and $H:(\C^N,p)\to (\C^N,p')$ a germ
of a holomorphic map sending $M$ into $M'$ whose Jacobian
does not vanish identically. Then for every positive
integer $\kappa$, there exists a germ of an algebraic
holomorphic map
$H^{\kappa}:(\C^N,p)\to (\C^N,p')$ which sends $M$ into
$M'$ and agrees with $H$ up to order $\kappa$.
\end{Thm}
One should  again note that the assumptions of Theorem
\ref{approx2} do not imply that the given holomorphic map
$H$ is itself algebraic. Theorem \ref{approx2}
immediately implies the following result concerning
biholomorphic and algebraic equivalence of generic
real-algebraic submanifolds.
\begin{Cor}\Label{alg}
Let $M,M'\subset \C^N$ be two real-algebraic generic
submanifolds of the same dimension. Assume that $M$ is
connected and of finite type at some point. Then for every
$p\in M$ and every $p'\in M'$, the germs $(M,p)$ and
$(M',p')$ are biholomorphically equivalent if and only if
they are algebraically equivalent.
\end{Cor}

In the case of real-algebraic hypersurfaces, we are able to drop the
finite type condition in Corollary \ref{alg}. In fact, we
shall prove the following.

\begin{Cor}\Label{hypalg}
Two germs of real-algebraic hypersurfaces in $\C^N$ are
biholomorphically equivalent if and only if
they are algebraically equivalent.
\end{Cor}

For  a positive integer $k$ and a point $p$ in $\C^N$,
denote  by
$G^k(\C^N,p)$ the jet group of order $k$ of $\C^N$ at $p$.
An element
$j(Z)$ of this group can be viewed as a $\C^N$-valued
polynomial in
$Z$ of degree at most $k$, fixing $p$, and with
nonvanishing Jacobian at $p$. The multiplication of two such
elements consist of composition of mappings with the
resulting polynomial truncated up to degree $k$ (see e.g.\ 
\cite{GG}).
If
$ (M,p)$ is a germ of a smooth  generic submanifold in
$\C^N$, we denote by ${\mathcal F}(M,p)$ the group of formal
invertible mappings $H:(\C^N,p)\to (\C^N,p)$ sending $M$
into itself. If $M$ is moreover assumed to be real-analytic,
then the subgroup of ${\mathcal F (M,p)}$ consisting
of those mappings which are convergent will be denoted by
$\Aut (M,p)$, the {\em stability group} of $M$ at
$p$. For any formal  map $H:(\C^N,p)\to
(\C^{N'},p')$, we define its jet $j_p^kH$ to
be its Taylor polynomial of degree
$k$ at $p$. If $N=N'$, $p=p'$ and $H$ is invertible, then
 $j_p^kH$ may be considered as an element of
$G^k(\C^N,p)$. The following corollary is a consequence of
Theorem \ref{fd} and Theorem
\ref{conv}.

\begin{Cor}\Label{inj}
Let $M\subset \C^N$ be a smooth generic submanifold with
$p\in M$. If $M$ is of finite type and
holomorphically nondegenerate at
$p$, then there exists a positive integer $K$ such that
the mapping $j_p^K:{\mathcal F}(M,p)\to G^K(\C^N,p)$ is
injective. If, in addition, $M$ is real-analytic, then
${\mathcal F}(M,p)=\Aut (M,p)$.
\end{Cor}

We shall now briefly mention previous work closely
related to the results in this paper. For the case of
Levi nondegenerate real-analytic hypersurfaces, finite
determination by their 2-jets and
convergence of formal invertible
maps were established in the seminal
paper of Chern-Moser \cite{CM}
(see also earlier work of Cartan \cite{Cartan} and
Tanaka \cite{ta}). The first and
third authors, jointly with Ebenfelt \cite{JAMS}
recently proved  the analogues of Theorems \ref{fd} and
\ref{conv} under the more restrictive condition that
$M'$ is essentially finite at $p'$, rather than
just holomorphically nondegenerate.  Earlier work by the
same authors on these topics appeared in \cite{CAG},
\cite{Asian}, \cite{MA}. The
second author of this paper established Theorem
1.2  (actually the more general version, Theorem 2.6
below) for the case of an invertible map $H$ between
real-analytic hypersurfaces \cite{Mirh}.
Theorem 2.6 was also proved
by the second author for invertible mappings
between generic real-analytic submanifolds of any
codimension under the additional assumption that one
of the manifolds is real-algebraic \cite{Mirg}. In
another direction, Ebenfelt \cite{e}  obtained  results
on finite determination  (not covered by Theorem
\ref{fd}) for smooth CR mappings between smooth
hypersurfaces. Lamel \cite{la} proved finite
determination and convergence results for certain
mappings between generic submanifolds of different
dimensions.  We also mention here that Theorem 1.2, for
the case of an invertible map  was
 claimed by Merker in a
preprint, http://xxx.lanl.gov/abs/math/9901027v1, but
the proof was incorrect (see Closing Remark in
\cite{JAMS}). 
While the present work was under
completion, Theorem 1.2 (in the form of Theorem 2.6)
was claimed a second time by Merker for invertible maps in the
preprint http://xxx.lanl.gov/abs/math/0005290v1, and
then again in a revision of that preprint in 
http://xxx.lanl.gov/abs/math/0005290v2. The authors
of the present paper have been unable to check the
proofs of the latter two preprints. 

It follows from \cite{CM} that if two germs of
real-analytic Levi nondegenerate hypersurfaces in
$\C^N$ are formally equivalent, then they are
biholomorphically equivalent.  On the other hand,
examples due to Moser and Webster \cite{MW} show that
there are pairs  of real-analytic
submanifolds which are formally
equivalent but are not biholomorphically equivalent. 
The first and third authors, in joint work with
Zaitsev \cite{BRZ} proved that, at ``general" points,
formal equivalence of real-analytic
submanifolds implies biholomorphic equivalence. 
Corollary \ref{fb} above establishes this result for points
not covered in previous published work. A related question for
real-algebraic submanifolds is the following, which has
been asked in
\cite{BAMS}:  If two germs of real-algebraic submanifolds are
biholomorphically equivalent, are they also
algebraically equivalent? It is shown in \cite{BRZ2} that
at ``general" points the answer is positive. 
Corollaries
\ref{alg} and
\ref{hypalg} above give further positive results for some
classes of submanifolds, including all hypersurfaces.
A related question is when a germ of a
holomorphic map sending one real-algebraic
submanifold into another is itself algebraic.  The
latter question has a long history.  We 
mention here the work of Webster \cite{Webster} for invertible
maps between Levi nondegenerate hypersurfaces, and, for
more recent work, we refer the reader to
\cite {Huang}, \cite{Acta}, \cite {Mia}, \cite{Z2}, and  \cite{CMS}.

Our approach in the proofs of the
results of this paper lies in the study
of the so-called ``reflection ideal"
associated to a triple $(M,M',H)$, where
$M$ and $M'$ are (germs of)
smooth generic submanifolds in $\C^
N$ and $\C^{N'}$ respectively, and
$H$ is a formal map sending $M$
into $M'$.  Such an ideal lies in the ring of formal
power series in $N+N'$ indeterminates. (See \S \ref{ideal} for precise
definitions.) If
the source generic submanifold $M$ is of
finite type, we establish finite
determination of reflection ideals
associated to formal mappings
(Theorem \ref{fdi} below) with no
nondegeneracy condition on the
target manifold $M'$. In fact, we
prove such a result in the more
general setting of formal manifolds.
When the generic submanifolds are
real-analytic and the source manifold $M$ is
of finite type, we prove (Theorem \ref{convideal} below)  that the
reflection ideal has a set of convergent generators.
If the generic submanifolds $M$ and $M'$ are
moreover real-algebraic, the map $H$ is convergent, and the
connected source manifold $M$ is of
finite type at some point, we prove
(Theorem \ref{algideal} below) that the
reflection ideal has a set of algebraic
generators.
An important ingredient for the
proofs of the above three theorems is
the use of iterated Segre mappings,
introduced in \cite{Acta} (see also
\cite{JAG}), which has already been
applied to various mapping
problems. Another important tool in
the proofs is Artin's approximation
theorem \cite{artin1} and an algebraic
version of the latter in \cite{artin2}. 

An outline of the organization of this paper is as
follows.  In
\S\ref{ideal} we state the more general results on
reflection ideals from which the theorems stated
above in this introduction will follow.  Sections \S \ref{jet}
to \S\ref{trick} are devoted to preliminaries needed for the
proofs of Theorems \ref {fdi},
\ref {convideal} and \ref{algideal}, which are given
in \S\S \ref{close} --\ref{2.7}.  Some remarks and open questions
are given in \S \ref{open}. 

\section{Manifold ideals and reflection
ideals}\Label{ideal}
For $x=(x_1,\ldots,x_k)\in \C^k$, we denote by $\C \dbl x
\dbr$ the ring of formal power series in $x$ and by
$\C\{x\}$ the subring of convergent ones. Moreover, we write
${\mathcal A}\{x\}\subset \C\{x\}$  for the subring of 
algebraic functions (also called Nash functions). If $R$ is
any of the three rings defined above and  
$I\subset R$ is an ideal generated by
$s_1(x),\ldots,s_d(x)$, we shall use the notation
$s(x)=(s_1(x),\ldots,s_d(x))$ and write
$I=(s(x))$. An ideal $I\subset R$ is called a
{\em manifold ideal} if it has a set of generators with
linearly independent differentials at the origin. Observe
that any two sets of such generators have the same number of
elements. This number is called the {\em
codimension} of
$I$. The following elementary fact, whose proof is left to
the reader, will be used implicitly throughout this paper. 
\begin{Lem}\Label{algebra}
Let $I\subset R$ be a manifold ideal of codimension $d$.
\begin{enumerate}
\item [(i)] Any set of $d$ elements of $I$ whose
differentials are linearly independent at
the origin generate
$I$.
\item [(ii)] From any set of generators
of
$I$, one may extract a subset of $d$ elements with
linearly independent differentials at the origin (which
generate
$I$ by {\rm (i)}).
\end{enumerate}
\end{Lem}

 If
$R$ is $\C\{x\}$ (respectively ${\mathcal A}\{x\}$)  and
$\{s_1(x),\ldots,s_d(x)\}$ is a set of generators of $I$ in
$R$, with $d$ the codimension of $I$, then the equations
$s_1(x)=\cdots=s_d(x)=0$ define a germ at 0 of a
complex-analytic (resp.\ complex-algebraic) submanifold
$\Sigma$ of codimension $d$. In general, we say that a
manifold ideal $I\subset \C\dbl x\dbr$ of codimension
$d$ defines a {\em formal manifold} $\Sigma \subset \C^k$
of dimension $k-d$ and write $I={\mathcal I}(\Sigma)$. (We
should point out that $\Sigma$ does not necessarily
correspond to a subset of $\C^k$ but we shall use the
notation $\Sigma
\subset
\C^k$ for motivation.) If $\Sigma \subset \C^k$ is a formal
manifold of dimension $l$, a parametrization of $\Sigma$ is
a formal mapping $(\C^l,0)\ni t\to v(t)\in (\C^k,0)$ such
that for any $h\in {\mathcal I}(\Sigma)$, $h\circ v=0$ and
$\rk
\partial v/\partial t(0)=l$.

If $I\subset \C \dbl x\dbr$ is an ideal and $F:(\C_x^k,0)\to
(\C^{k'}_{x'},0)$ is a formal map, then  the {\em
pushforward} $F_*(I)$ of $I$ is defined to be the ideal in
$\C \dbl x'\dbr$, 
$x'\in \C^{k'}$, 
\begin{equation}\Label{push}
F_*(I):=\{h\in \C \dbl x'\dbr:h\circ F\in
I\}.
\end{equation}
If $\Sigma \subset \C^k$ and $\Sigma'\subset
\C^{k'}$ are formal manifolds with $I={\mathcal I}(\Sigma)\subset
\C\dbl x\dbr$ and $I'={\mathcal I}(\Sigma')\subset
\C\dbl x'\dbr$, then we say that $F$ sends $\Sigma$ into 
$\Sigma'$ and write $F(\Sigma)\subset \Sigma'$ if
$I'\subset F_*(I)$. 

For a formal map $F:(\C^k_x,0)\to (\C^{k'}_{x'},0)$, we
denote by $\Rk F$ the rank of the
Jacobian  matrix $\d F /\d x$ regarded as a
$\C \dbl x\dbr$-linear mapping
$(\C \dbl x\dbr)^k\to (\C \dbl x\dbr)^{k'}$.
Hence
$\Rk F$ is the largest integer $r$ such that there
is an
$r\times r$ minor of the matrix $\d F/\d x$ which is
not  0 as a formal power series in $x$. Note that if $F$ is
convergent, then $\Rk F$ is the generic rank of the map $F$.

\begin{Def}\Label{nono}
{\rm Let $\Sigma\subset \C^k$ and
$\Sigma'\subset \C^{k'}$ be two formal manifolds of
dimension $l,l'$ respectively and
$F:(\C^k,0)\to (\C^{k'},0)$ a formal map sending $\Sigma$
to $\Sigma'$. Then $F$ is said to be {\it
$(\Sigma,\Sigma')$-nondegenerate} if $\Rk F\circ
v=l'$ for some (and hence for all) parametrization $v$ of
$\Sigma$.}
\end{Def}
A formal vector field $V$ in $\C^k$ is a $\C$-linear
derivation of $\C\dbl x\dbr$ and hence is given by 
$$V=\sum_{j=1}^ku_j(x)\frac{\partial}{\partial x_j},\quad
u_j(x)\in \C \dbl x\dbr,\ j=1,\ldots,k.$$
The vector field $V$ is called tangent to a formal manifold
$\Sigma \subset \C^k$ or, equivalently, to its ideal $\I
(\Sigma)$ if and only if $V(f)$ belongs to $\I (\Sigma)$ for
every $f\in \I (\Sigma)$.

\begin{Def}\Label{pear}
{\rm 
An ideal $I\subset \C\dbl x\dbr$ is said to be
{\em convergent} (resp.\ {\em algebraic}) if
$I$ has a set of convergent (resp.\ algebraic) generators.}
\end{Def}

 For
$(Z,\zeta)\in
\C^N
\times
\C^N$, we define the involution $\sigma : \C\dbl
Z,\zeta\dbr \to
\C\dbl Z,\zeta\dbr$ by $\sigma
(f)(Z,\zeta):=\bar{f}(\zeta,Z)$, where $\bar{f}$ is the
formal power series obtained from $f$ by taking complex
conjugates of the coefficients. An ideal
${\mathcal J}\subset \C\dbl Z,\zeta\dbr$ is called {\em
real} if $\sigma (f)\in {\mathcal J}$ for every $f\in
{\mathcal J}$. Since $\sigma$ is also an involution when
restricted to $\C\{Z,\zeta\}$ or ${\mathcal A}\{Z,\zeta\}$,
a similar definition applies for ideals in these rings.
A formal manifold $\scrM \subset \C^N\times \C^N$ is called
{\em real} if its ideal $\I (\scrM)$ is real.
A  formal real manifold $\scrM \subset \C^N\times \C^N$ of
codimension $d$ is called {\em generic} if for some (and
hence for any) vector of $d$ generators
$\rho (Z,\zeta)=(\rho_1(Z,\zeta),\ldots,\rho_d(Z,\zeta))$ of
$\I (\scrM)$, the rank of the $d\times N$ matrix $\partial
\rho/\partial Z(0)$ is $d$. To motivate this definition, let
$M\subset
\C^N$ be a smooth generic submanifold of codimension $d$
through the origin with smooth local defining functions
$r(Z,\bar{Z})=(r_1(Z,\bar{Z}),\ldots,r_d(Z,\bar{Z}))$ whose
 Taylor expansions at zero are $\rho
(Z,\bar{Z})=(\rho_1(Z,\bar{Z}),\ldots,\rho_d(Z,\bar{Z}))$.
Observe that the
$d$ vector-valued formal power series $\rho (Z,\zeta)$
generate a real manifold ideal in $\CZZ$ whose formal
manifold $\scrM
\subset
\C^N\times \C^N$ is generic. If, furthermore, $M$ is
real-analytic, then $\scrM \subset \C^N\times \C^N$ is a
germ at 0 of a complex submanifold of codimension $d$,
usually referred to as the {\em complexification} of $M$.

For a formal generic manifold ${\scrM}\subset \C^N\times
\C^N$, we define a manifold ideal ${\mathcal
I}_0(\scrM)\subset \C \dbl Z \dbr$ as the ideal generated by
the $h(Z,0)$ for all $h\in {\mathcal I}(\scrM)$. 
The formal manifold $S_0(\scrM)\subset
\C^N$ associated to this ideal is called the {\em
formal Segre variety} of
${\scrM}$ at 0. Observe that when $\scrM$ is the
complexification of a real-analytic generic submanifold
$M\subset \C^{N}$ (through 0), then $S_0(\scrM)$ is the
usual Segre variety of $M$ at 0. 

For a formal map
$H:(\C_Z^N,0)\to (\C_{Z'}^{N'},0)$, we define its {\em
complexification} ${\mathcal H}:(\C_{Z}^N\times
\C^{N}_{\zeta},0)\to (\C^{N'}_{Z'}\times
\C^{N'}_{\zeta'},0)$ to be the formal map given by 
\begin{equation}\Label{H}
{\mathcal H}(Z,\zeta):=
(H(Z),\bar{H}(\zeta)).
\end{equation}
In what follows, given 
${\scrM}\subset
\C^N\times
\C^{N}$ and
$\scrM'\subset \C^{N'}\times \C^{N'}$ two formal generic
manifolds, we will consider formal maps
$H:(\C^N,0)\to (\C^{N'},0)$ such that their
complexifications ${\mathcal H}$, as defined by (\ref{H}),
send
${\scrM}$ into $\scrM'$.  It is
easy to check that if
$H$ is such a mapping, then $H$ sends the formal Segre
variety $S_0(\scrM)$ into the formal Segre variety
$S_0(\scrM')$. 
\begin{Def}\Label{apple} {\rm Let $\scrM\subset
\C^N\times \C^N$ and $\scrM'\subset
\C^{N'}\times \C^{N'}$ be two formal generic
manifolds and $H:(\C^N,0)\to (\C^{N'},0)$ a formal map
such that its complexification
${\mathcal H}$ maps $\scrM$ into $\scrM'$. 
The map $H$ is called {\em not totally degenerate} if $H$ is
$(S_0(\scrM),S_0(\scrM'))$-nondegenerate as defined in
 Definition \ref{nono}.}
\end{Def}

A formal (1,0)-vector field $X$ in $\C_Z^N \times
\C_{\zeta}^N$ is given by
\begin{equation}\Label{1,0}
X=\sum_{j=1}^Na_j(Z,\zeta)\frac{\partial}{\partial
Z_j},\quad a_j(Z,\zeta)\in \CZZ,\ j=1,\ldots,N.
\end{equation}
Similarly, a (0,1)-vector field $Y$ in $\C_{Z}^N \times
\C_{\zeta}^N$ is given by
\begin{equation}\Label{0,1}
Y=\sum_{j=1}^Nb_j(Z,\zeta)\frac{\partial}{\partial
\zeta_j},\quad b_j(Z,\zeta)\in \CZZ,\ j=1,\ldots,N.
\end{equation}
For a formal generic manifold $\scrM \subset \C^N\times
\C^N$ of codimension $d$, we denote by $\frak g_{\scrM}$ 
the Lie algebra generated by the formal (1,0) and (0,1)
vector fields tangent to $\scrM$. The formal generic
manifold
$\scrM$ is said to be {\em of finite type} if the dimension
of
$\frak g_{\scrM}(0)$ over $\C$ is $2N-d$, where $\frak
g_{\scrM}(0)$ is the vector space obtained by evaluating
the vector fields in $\frak g_{\scrM}$ at the origin of
$\C^{2N}$.

Let $H:(\C^N,0)\to (\C^{N'},0)$ be a formal mapping. For an
ideal $J\subset \C \dbl Z',\zeta'\dbr$, $(Z',\zeta')\in
\C^{N'}\times \C^{N'}$, 
we define $J^H\subset \C \dbl Z,\zeta' \dbr$ to be the
ideal generated by the $h(H(Z),\zeta')$ for all $h\in
J$ i.e. 
\begin{equation}\Label{ballot}
J^H:=(h(H(Z),\zeta'):h\in
J)\subset
\C
\dbl Z,\zeta'
\dbr.
\end{equation}
Note that if $J$ is generated by
$s(Z',\zeta')=(s_{1}(Z',\zeta'),\ldots,s_{m}(Z',\zeta'))$
in $\C \dbl Z',\zeta'\dbr$, then $J^H$ is generated by the
components of $s(H(Z),\zeta')$ in $\C \dbl Z,\zeta'\dbr$. If
$\scrM'\subset
\C_{Z'}^{N'}\times
\C_{\zeta'}^{N'}$ is a formal generic submanifold of
codimension $d'$, we write for simplicity of notation
\begin{equation}\Label{crap}
\I^H:=\I (\scrM')^H\subset \C \dbl Z,\zeta'
\dbr,
\end{equation}
where we have used the notation given in (\ref{ballot}).
	It is easy to see that $\I^H$ is a manifold ideal of
codimension
$d'$ in $\C\dbl Z,\zeta'\dbr$. If ${\scrM'}$
and $H$ are as above, then we
refer to the ideal ${\mathcal I}^H$ as the
{\em reflection ideal of $H$} (relative to ${\scrM'}$). If
$(M',0)$ is a germ of a real-analytic generic submanifold
of $\C^{N'}$ and $H:(\C^N,0)\to (\C^{N'},0)$ is a formal
map, we again define ${\mathcal I}^H$ by (\ref{crap}), where
${\scrM'}$ is the complexification of $M'$.

 Our
first result in this section establishes finite
determination of reflection ideals for formal mappings $H$
such that their complexifications
${\mathcal H}$ defined in (\ref{H}) send a formal generic
manifold
$\scrM$ into $\scrM'$. Note that in Theorem \ref{fdi}, no
nondegeneracy condition is imposed on the formal manifold
$\scrM'$.

\begin{Thm}\Label{fdi}
Let $\scrM \subset \C^N\times \C^N$ and $\scrM'\subset
\C^{N'}\times \C^{N'}$ be formal generic manifolds with
$\scrM$ of finite type. Let $H^0:(\C^N,0)\to (\C^{N'},0)$
be a formal map such that its complexification ${\mathcal
H}^0$ sends
$\scrM$ into $\scrM'$. Assume furthermore that $H^0$ is not
totally degenerate as in Definition {\rm \ref{apple}}.
Then there exists a positive integer
$K_0$ such that if $H:(\C^N,0)\to (\C^{N'},0)$ is a formal
map with ${\mathcal H}(\scrM)\subset \scrM'$ and
$j_0^{K_0}H=j_0^{K_0}H^0$, it follows that the
corresponding reflection ideals defined by {\rm
(\ref{crap})} are the same i.e. 
\begin{equation}\Label{eqideals}
{\mathcal I}^H={\mathcal I}^{H^0}.
\end{equation}
\end{Thm}

If $(M,0)$ and $(M',0)$ are germs of real-analytic generic
submanifolds in $\C^N$ and $\C^{N'}$ respectively and
$H:(\C^N,0)\to (\C^{N'},0)$ is a formal mapping sending $M$
into $M'$ as defined in \S \ref{int}, then its
complexification
${\mathcal H}$ sends
${\scrM}$ into $\scrM'$, where ${\scrM}$ and $\scrM'$ are
the complexifications of $M$ and $M'$ respectively. The
second main result of this section establishes convergence
of reflection ideals for formal mappings between
real-analytic generic submanifolds, with no 
nondegeneracy condition imposed on the target manifold $M'$.

\begin{Thm}\Label{convideal}
Let $(M,0)$ and $(M',0)$ be germs of real-analytic generic
submanifolds in $\C^N$ and $\C^{N'}$ respectively and
$H:(\C^N,0)\to (\C^{N'},0)$ a formal mapping sending $M$
into $M'$. Assume that $M$ is of finite type at $0$ and $H$
is not totally degenerate. Then  the reflection
ideal ${\mathcal I}^H$, as defined by {\rm
(\ref{crap})}, is convergent.
\end{Thm}

The
last result of this section establishes algebraicity
of reflection ideals for local holomorphic mappings between
real-algebraic generic submanifolds, with no 
nondegeneracy condition imposed on the target manifold $M'$.
\begin{Thm}\Label{algideal}
Let $M,M'\subset \C^N$ be real-algebraic generic
submanifolds of codimension $d$ through the origin and
$H:(\C^N,0)\to (\C^N,0)$ be a germ of a holomorphic map
sending $M$ into $M'$. Assume that the Jacobian of $H$ does
not vanish identically and that there is no germ of a
nonconstant holomorphic function
$h:(\C^N,0)\to
\C$ with
$h(M)\subset
\bR$.  Then the reflection
ideal ${\mathcal I}^H$, as defined by {\rm
(\ref{crap})}, is algebraic.
\end{Thm}
In view of Proposition \ref{both} (iii) below, Theorem
\ref{algideal} in the case where $M$ and $M'$
are real-algebraic hypersurfaces in $\C^N$, is contained in
\cite{Mia}.
\begin{Rem}\Label{tex}
{\rm Even if all the assumptions of Theorem \ref{convideal}
are satisfied, the fact that the reflection ideal
${\mathcal I}^H$ is convergent does not imply that the
formal map $H$ is convergent. For example, let $M=M'$ be the
real-algebraic hypersurface of finite type through the
origin in
$\C^3$ given by 
$${\rm Im}\ Z_3=|Z_1Z_2|^2.$$
For any nonconvergent
formal power series  $h(Z)=h(Z_1,Z_2,Z_3)$ vanishing at the
origin, let $H:(\C^3,0)\to (\C^3,0)$ be the formal 
invertible map  given by
$$H(Z_1,Z_2,Z_3)=
(Z_1e^{h(Z)},Z_2e^{-h(Z)},Z_3).$$
In this example, the formal map $H$ sends $M$ into itself
and is not convergent, but one can easily check that its
reflection ideal
${\mathcal I}^{H}$ is convergent. (This fact  also follows
from Theorem
\ref{convideal}.) Similar considerations can be made in the
algebraic case relative to Theorem \ref{algideal}.
Proposition \ref{convhol} below gives an additional
condition on $M'$ which guarantees that the convergence of
${\mathcal I}^H$ implies that $H$ is convergent. }
\end{Rem}
The following proposition, which justifies the notion
of convergent reflection ideals introduced here, will be
used for the proofs of Theorems \ref{approx1} and
\ref{approx2}.
\begin{Pro}\Label{link}
Let $(M',0)$ be a germ of a generic real-analytic (resp.\ 
real-algebraic) submanifold of codimension
$d'$ in $\C^{N'}$ and $H:(\C^N,0)\to (\C^{N'},0)$ a
formal map. Then the reflection ideal ${\mathcal I}^H$ is
convergent (resp.\  algebraic) if and only if there exists a
convergent (resp.\  algebraic) map $\check H:(\C^N,0)\to
(\C^{N'},0)$ such that ${\mathcal I}^H={\mathcal I}^{\check
H}$. More precisely, if
${\mathcal I}^H$ is convergent (resp.\ algebraic), then for
any positive integer $\kappa$, there exists a
convergent (resp.\ algebraic) map $H^{\kappa}:(\C^N,0)\to
(\C^{N'},0)$ agreeing up to order $\kappa$ with $H$ such
that
${\mathcal I}^H={\mathcal I}^{H^{\kappa}}$.
\end{Pro}

If ${\scrM}\subset \C^N\times \C^N$ is a formal generic
manifold, we say that ${\mathcal M}$ is {\em holomorphically
nondegenerate} if there is no nontrivial (1,0) vector field
of the form (\ref{1,0}) tangent to ${\scrM}$ with
coefficients $a_j(Z,\zeta)=a_j(Z)$ independent of $\zeta$
for $j=1,\ldots,N$. Note that if $(M,0)$ is a germ of
a smooth generic submanifold in $\C^N$, then $M$ is
holomorphically nondegenerate in the sense defined in \S
\ref{int} if and only if its associated formal generic
manifold $\scrM$ is holomorphically nondegenerate as
defined here. If $\scrM$ is the
complexification of a germ of a real-analytic
generic submanifold $(M,0)$ in $\C^N$, then $\scrM$ is
holomorphically nondegenerate as defined here if and only if
there is no germ of a nontrivial (1,0) vector field
of the form (\ref{1,0}) tangent to ${\scrM}$ with
convergent coefficients $a_j(Z,\zeta)=a_j(Z)$ independent of
$\zeta$ for $j=1,\ldots,N$ (see e.g.\ \cite{BER}).

Theorem \ref{fdi} will be used in conjunction with the
following finite determination result to prove Theorem
\ref{fd}.

\begin{Pro}\Label{fdhol}
Let ${\scrM'}\subset \C^{N'}\times \C^{N'}$ be a
holomorphically nondegenerate formal generic manifold and
$H^0:(\C^N,0)\to (\C^{N'},0)$ a formal map with $\Rk
H^0=N'$. Then there exists a positive integer $K$ such that
if
$H:(\C^N,0)\to (\C^{N'},0)$ is a formal map with
$j_0^KH=j_0^KH^0$ and
${\mathcal I}^H={\mathcal I}^{H^0}$ (as defined in {\rm
(\ref{crap})}), it follows that
$H=H^0$.
\end{Pro}

In the case of a real-analytic generic submanifold and
holomorphic mappings, we have the following
geometric interpretation of the equality (\ref{eqideals}) of
reflection ideals. In view of Proposition \ref{geom}
below, Theorem \ref{fdi} can then be seen as a finite
determination result for Segre varieties. 

\begin{Pro}\Label{geom}
Let $(M',0)$ be a germ of a real-analytic generic
submanifold in $\C^{N'}$ with real-analytic
local defining functions
$r'(Z',\bar{Z}')$. Assume that
$H,H^0:(\C^N,0)\to (\C^{N'},0)$ are germs of holomorphic
mappings. Then the following two conditions are equivalent:
\begin{enumerate}
\item [(i)] ${\mathcal I}^H={\mathcal I}^{H^0}$, where
the ideals ${\mathcal I}^{H}$ and ${\mathcal I}^{H^0}$ are
defined by {\rm (\ref{crap})}.
\item [(ii)] For $Z$ near the origin, the Segre varieties
of $M'$ relative to the points $H(Z)$ and $H^0(Z)$ are the
same. More precisely, there exists open neighborhoods of
$0$,
$U$ and
$U'$ in
$\C^N$ and
$\C^{N'}$ respectively such that for all $Z\in U$,
\begin{equation}
S_{H(Z)}=S_{H^0(Z)},
\end{equation}
where
$S_{H(Z)}=\{Z'\in U':r'(Z',\overline{H(Z)})=0\}$, with
a similar definition for $S_{H^0(Z)}$.
\end{enumerate}
\end{Pro}
Proposition \ref{geom} will not be used in
the remainder of the paper and its proof is left to the
reader. The last result of this section connects the
convergence of the reflection ideal ${\mathcal I}^H$ to the
convergence of the mapping $H$.
\begin{Pro}\Label{convhol}
Let $(M',0)$ be a germ of a generic real-analytic
holomorphically nondegenerate submanifold of codimension
$d'$ in
$\C^{N'}$. If
$H:(\C^N,0)\to (\C^{N'},0)$ is a formal map with $\Rk
H=N'$ such that its reflection ideal ${\mathcal
I}^{H}$, as defined  by {\rm (\ref{crap})}, is convergent,
then
$H$ is convergent.
\end{Pro}

\begin{Rem}\Label{UTC}
{\rm We should point out that a statement similar to
Proposition \ref{convhol} holds in the algebraic case.
Indeed, if $(M',0)$ is a germ of a generic real-algebraic
holomorphically nondegenerate submanifold of codimension
$d'$ in
$\C^{N'}$ and if
$H:(\C^N,0)\to (\C^{N'},0)$ is a formal map with $\Rk
H=N'$ such that its reflection ideal ${\mathcal
I}^{H}$ is algebraic,
then
$H$ is algebraic. This fact will not be used in this paper.}
\end{Rem}
The proofs of Proposition \ref{link}, \ref{fdhol} and
\ref{convhol} will be given in \S \ref{tenure}.

\section{Further results on finite determination, 
convergence, and approximation of mappings}\Label{blou}

The following finite determination result,  which is a
generalization of Theorem \ref{fd}, will be a consequence of
Theorem
\ref{fdi} and Proposition \ref{fdhol}.

\begin{Thm}\Label{plus1}
Let $\scrM \subset \C^N\times \C^N$ and $\scrM'\subset
\C^{N'}\times \C^{N'}$ be formal generic manifolds with
$\scrM$ of finite type and $\scrM'$ holomorphically
nondegenerate. Let
$H^0:(\C^N,0)\to (\C^{N'},0)$ be a formal map such that its
complexification
${\mathcal H}^0$ sends
$\scrM$ into $\scrM'$. Assume furthermore that $H^0$ is not
totally degenerate as in Definition {\rm \ref{apple}} and
that $\Rk H^0=N'$. Then there exists a positive integer
$K$ such that if $H:(\C^N,0)\to (\C^{N'},0)$ is a formal
map with ${\mathcal H}(\scrM)\subset \scrM'$ and
$j_0^{K}H=j_0^{K}H^0$, it follows that $H=H^0$.
\end{Thm}

Similarly, the following convergence result,  which
is a generalization of Theorem \ref{conv}, will be a
consequence of Theorem
\ref{convideal} and Proposition \ref{convhol}.

\begin{Thm}\Label{plus2}
Let $(M,0)$ and $(M',0)$ be germs of real-analytic generic
submanifolds in $\C^N$ and $\C^{N'}$ respectively and
$H:(\C^N,0)\to (\C^{N'},0)$ a formal mapping sending $M$
into $M'$. Assume that $M$ is of finite type at $0$ and
that $M'$ is holomorphically nondegenerate at $0$. If
$H$ is not totally degenerate and $\Rk H=N'$, then  
$H$ is convergent.
\end{Thm}

\begin{Rem}
{\rm We should point out that the assumptions of
Theorem \ref{plus2} are less restrictive than those of
Theorem \ref{conv}, even in  the
case where $M$ and $M'$ are real-analytic hypersurfaces in
the same space $\C^N$. (The same can also be said about
Theorems \ref{plus1} and \ref{fd}.) For instance,
given a nontrivial convergent power series
$h:(\C,0)\to (\C,0)$, consider the  following hypersurfaces
in $\C^3$:
\begin{equation}
\begin{aligned}
M:=& \{Z\in \C^3: \im
Z_3=|Z_1^2Z_2|^2+|h(Z_1)|^2
\},\\
M':=& \{Z'\in \C^3: \im
Z'_3=|Z'_1Z'_2|^2+|h(Z'_1)|^2\}.
\end{aligned}
\end{equation}
Observe that the convergent mapping $(\C^3,0)\ni Z\mapsto
H(Z):=(Z_1,Z_1Z_2,Z_3)\in (\C^3,0)$ sends $M$ into $M'$.
Moreover, $M$ and $M'$ are of finite
type
and holomorphically nondegenerate at the origin. Note also
that $H$ is not totally degenerate and $\Rk H=3$ but $H$ is
not finite. We should point out that the convergence of
formal mappings between
$M$ and
$M'$ satisfying the latter conditions follows from Theorem
\ref{plus2}, but does not follow  from Theorem
\ref{conv} nor from previously known  results. (Indeed,
since
$M$ and
$M'$ are not essentially finite at the origin, the result
in \cite{JAMS} does not apply, nor does the one in
\cite{Mirg} if the function $h$ is chosen not be
algebraic.)}
\end{Rem}

The following approximation result generalizes Theorem
\ref{approx1}.

\begin{Thm}\Label{plus3}
Let $(M,0)$ and $(M',0)$ be two germs of real-analytic
generic submanifolds in $\C^N$ and $\C^{N'}$ respectively,
with
$M$ of finite type at $0$. If $H:(\C^N,0)\to (\C^{N'},0)$
is a not totally degenerate formal map sending $M$ into
$M'$ and if
$\kappa$ is a positive integer, then there exists a
convergent map
$H^{\kappa}:(\C^N,0)\to (\C^{N'},0)$ which sends
$M$ into $M'$ and agrees with $H$ up to order
$\kappa$.
\end{Thm}

\section{Ideals in jet spaces}\Label{jet}

Given nonnegative integers $l,k,r$, with $k,r\geq 1$, we
denote by
$J^l_0(\C^k,\C^r)$ the jet space at the origin of order
$l$ of holomorphic mappings from $\C^k$ to $\C^r$. An
element $j$ of $J^l_0(\C^k,\C^r)$ can be written as a
polynomial mapping 
\begin{equation}\Label{coord}
j(X)=\sum_{\alpha \in
\Nn ^k,\ 0\leq |\alpha|\leq
l}\frac{\Lambda_{\alpha}}{\alpha!}X^{\alpha},\quad
\Lambda_{\alpha} \in \C^{r}.
\end{equation}
We think of the coefficients
$\Lambda:=(\Lambda_{\alpha})_{0\leq |\alpha|\leq l}$,
$\Lambda_{\alpha}\in \C^r$, as linear  coordinates in the
finite dimensional vector space
$J^l_0(\C^k,\C^r)$ and we identify $j$ with $\Lambda$. 
We write $\Lambda_{\alpha}=(\Lambda_{\alpha,i})_{1\leq i
\leq r}$ for any $\alpha \in \N^k$, $|\alpha|\leq l$. We
also use the splitting
\begin{equation}\Label{Lambda}
\Lambda=(\Lambda_0,\hat \Lambda),\quad
\hat \Lambda=(\Lambda_{\alpha})_{1\leq |\alpha|\leq l}.
\end{equation} 
 Using
the coordinates $\Lambda$, we identify $J^l_0(\C^k,\C^r)$
with
$\C^m_{\Lambda}$ where
$m={\rm dim}_{\C}J^l_0(\C^k,\C^r)$.  For a formal map
$F:(\C^k,0)\to (\C^r,0)$, we write
$j_x^lF$ and $\hat \jmath^l_xF$ for the vectors of formal
series 
\begin{equation}\Label{number}
j_x^lF:=(\partial^{\nu}F(x))_{0\leq
|\nu|\leq l},\quad
\hat \jmath^l_xF:=(\partial^{\nu}F(x))_{1\leq |\nu|\leq l}.
\end{equation}
Here, for $\nu \in
\N^k$, 
$\partial^{\nu}F(x)\in (\C\dbl x\dbr)^r$ and $x\in \C^k$.
If $s$ is another positive integer and
$\eta:(J^l_0(\C^k,\C^r),0)\to (J^l_0(\C^k,\C^s),0)$ is a
formal map, we take coordinates
$\Lambda=(\Lambda_{\nu})_{|\nu|\leq l}$ and
$\Lambda'=(\Lambda'_{\nu})_{|\nu|\leq l}$ for
$J^l_0(\C^k,\C^r)$ and $J^l_0(\C^k,\C^s)$ respectively.
Here, $\Lambda_{\nu}\in \C^r$, $\Lambda'_{\nu}\in \C^s$. We
then write $\eta_{\nu}=\Lambda'_{\nu}\circ \eta$; that is, 
the map $\eta$ is given by $\Lambda'=\eta (\Lambda)$.
Hence, for $\nu \in \N^k$, $|\nu|\leq l$, $\eta_{\nu}$ is
the $\nu$-th component of $\eta$ i.e.
\begin{equation}\Label{crappy}
\eta_{\nu}:(J^l_0(\C^k,\C^r),0)\to
(\C^s,0),\quad \eta=(\eta_{\nu})_{|\nu|\leq l}.
\end{equation}

If $R$ is a ring and $T\in \C^q$, as usual we denote by 
$R[T]$ the ring of polynomials in $T$ with coefficients in
$R$. If $\Lambda=(\Lambda_0,\hat \Lambda)$ are coordinates
in $J_0^l(\C^k,\C^r)$ as in (\ref{coord}) and
(\ref{Lambda}), the subring 
$\C \dbl \Lambda_0\dbr[\hat \Lambda]:=(\C \dbl
\Lambda_0\dbr)[\hat \Lambda]$ of the ring $\C \dbl
\Lambda\dbr$ will play a crucial role in the rest of this
paper. For instance, if $u\in \C
\dbl
\Lambda_0\dbr [\hat \Lambda]$ and $F:(\C^k,0)\to (\C^r,0)$
is a formal map, then $u(j_x^lF)$ is a well-defined formal
power series in $\C \dbl x\dbr$ (while for a general $u\in
\C
\dbl
\Lambda\dbr$, one cannot define $u(j_x^lF)$). We have the
following uniqueness result.

\begin{Lem}\Label{lema}
If $\Lambda=(\Lambda_0,\hat \Lambda)$ are coordinates
in $J_0^l(\C^k,\C^r)$ as in {\rm (\ref{coord})} and {\rm
(\ref{Lambda})}  and if
$u\in
\C
\dbl
\Lambda_0\dbr [\hat \Lambda]$ is a formal power series
satisfying
\begin{equation}\Label{dimple}
u(j_x^lF)=0,\ {\rm in }\ \C \dbl
x\dbr,
\end{equation}  for any formal power series mapping
$F:(\C^k,0)\to (\C^r,0)$, then
$u=0$ $($in $\C \dbl \Lambda_0\dbr [\hat \Lambda]$$)$.
\end{Lem}

\begin{proof} We shall define a polynomial map 
\begin{equation}\Label{org}
\varphi :(J_0^l(\C^k,\C^r)\times \C^k,0)\to
(J_0^l(\C^k,\C^r),0)
\end{equation}
as follows. If $A=(A_{\nu})_{|\nu|\leq l}$, $\nu \in \N^k$,
$A_{\nu}\in \C^r$, are coordinates on the source jet space
$J_0^l(\C^k,\C^r)$ as in (\ref{coord}),
$x=(x_1,\ldots,x_k)\in
\C^k$ and $\Lambda=(\Lambda_{\alpha})_{|\alpha|\leq l}$ are
coordinates on the target jet space $J_0^l(\C^k,\C^r)$,
then $\varphi$ is defined by
\begin{equation}\Label{gold}
\Lambda=\varphi
(A,x):=\left(\partial_x^{\alpha}
\big(x_1\sum_{0\leq |\nu|\leq
l}A_{\nu}x^{\nu}\big)
\right)_{0\leq |\alpha|\leq l}.
\end{equation}
We claim that the generic rank of $\varphi$, $\Rk \varphi$,
is equal to $m$, the dimension of $J_0^l(\C^k,\C^r)$ over
$\C$. For this, let 
$$\tilde\varphi (A,x):=\left(x_1
\partial_x^{\alpha}\big(
\sum_{0\leq |\nu|\leq l}A_{\nu}x^{\nu}\big)\right)_{0\leq
|\alpha|\leq l}.$$
First note that $\Rk
\varphi$ is greater or equal to the generic rank (in $x$)
of the
$m\times m$ matrix $\Delta (x):=\displaystyle \frac{\partial
\varphi}{\partial A}(A,x)$. Moreover, it is
not difficult to see that the generic rank of the
matrix $\Delta (x)$ is the same as that of the matrix $\tilde
\Delta (x):=\displaystyle 
\frac{\partial \tilde \varphi
}{\partial A}(A,x)$. 
Since for $x_1\not =0$, the  rank of $\tilde \Delta(x)$
is clearly $m$, it follows that $\Rk \varphi =m$. This
proves the claim.

Given a formal power series $v\in \C\dbl \Lambda\dbr$,
since $\varphi (0)=0$, we can consider the composition
$(v\circ \varphi )(A,x)$ as a formal power series in $\C
\dbl A,x\dbr$. We write 
\begin{equation}\Label{poi}
(v\circ \varphi )(A,x)=\sum_{\beta \in
\Nn^k}v_{\beta}(A)x^{\beta},\quad v_{\beta}\in \C \dbl
A\dbr.
\end{equation}
Observe that if $v$ is in the subring $\C\dbl
\Lambda_0\dbr[\hat\Lambda]\subset \C \dbl \Lambda\dbr$ then
for each $\beta\in \N^k$, $v_{\beta}$ is a polynomial in
$A$, i.e.\ $(v\circ \varphi)(A,x)\in \C [A]\dbl x\dbr$. Let
$u$ be as in Lemma \ref{lema} satisfying (\ref{dimple}).
For any vector
$a=(a_{\nu})_{|\nu|\leq l}\in
\C^m\cong J_0^l(\C^k,\C^r)$,  by (\ref{dimple}) with 
$F(x)=x_1\sum_{0\le |\nu|\leq l}a_{\nu}x^{\nu}$, we obtain
\begin{equation}\Label{eq1}
u(j_x^lF)=u (\varphi (a,x))=\sum_{\beta \in
\Nn^k}u_{\beta}(a)x^{\beta}=0,\ {\rm in}\
\C
\dbl x\dbr.
\end{equation}
As a consequence, we have $u_{\beta}(a)=$ for any $\beta
\in \N^k$ and any vector $a$ in $\C^m$. Since $u_{\beta}$
is a polynomial, it follows that $u_{\beta}\equiv 0$ and
hence the formal power series $(u\circ \varphi)(A,x)$ is
zero in $\C [A]\dbl x\dbr\subset \C \dbl A,x\dbr$. To
conclude that
$u$ is identically zero, by e.g.\ Proposition 5.3.5 of
\cite{BER}, it suffices to use the fact that $\Rk
\varphi=m$. This completes the proof of Lemma \ref{lema}.
\end{proof}

\begin{Pro}\Label{gener} Let $l,r,s$ be nonnegative
integers with  $r,s\geq 1$, and let
$\phi:(\C^r,0)\to (\C^s,0)$ be a formal map. Then there
exists  a unique  formal map
\begin{equation}\Label{phil}
\phi^{(l)}:(J_0^l(\C^k,\C^r),0)\to
(J_0^l(\C^k,\C^s),0)
\end{equation} 
whose components are in $\C \dbl \Lambda_0\dbr [\hat
\Lambda]$, with $\Lambda =(\Lambda_0,\hat \Lambda)$ 
the coordinates of $J^l_0(\C^k,\C^r)$ introduced in {\rm
(\ref{coord})} and {\rm (\ref{Lambda})}, such that for any formal map
$F:(\C^k,0)\to (\C^r,0)$
\begin{equation}\Label{jetid}
j_x^l (\phi\circ F)=\phi^{(l)}(j_x^lF).
\end{equation}
Moreover, if we write
$\phi^{(l)}(\Lambda)=(\phi^{(l)}_{\nu}(\Lambda))_{\nu \in
\Nnn^k,|\nu|\leq l}$, then for each $\nu$,
$\phi^{(l)}_{\nu}(\Lambda)$ depends only on
$(\Lambda_{\alpha})_{\alpha \leq \nu}$.
Finally, if $r=s$
and $\phi:(\C^r,0)\to (\C^r,0)$ is invertible, then so is
$\phi^{(l)}:(J_0^l(\C^k,\C^r),0)\to (J_0^l(\C^k,\C^r),0)$
and
$(\phi^{(l)})^{-1}=(\phi^{-1})^{(l)}$.
\end{Pro}
\begin{proof} The existence of the map $\phi^{(l)}$ and its
properties follow easily from the chain rule. The
uniqueness of such a map is a consequence of Lemma
\ref{lema}. The proof of the last statement of the
proposition is straightforward and left to the reader.
\end{proof}

\begin{Rem}\Label{parameter}
{\rm Let $\phi$ and $\phi^{(l)}$ be as in Proposition
{\rm \ref{gener}}. It follows from {\rm (\ref{jetid})} and
the other properties of $\phi^{(l)}$ that for any
formal map $G:(\C_x^k\times
\C_t^q,0)\to (\C_y^r,0)$, we have the equality of
vector valued formal power series in
$\C\dbl x,t\dbr$
\begin{equation}\Label{eqparam}
j_x^l(\phi (G(x,t)))=\phi^{(l)}(j^l_xG(x,t)).
\end{equation}
Here, as in (\ref{number}),
$j^l_xG(x,t)=(\partial_x^{\nu}G(x,t))_{|\nu|\leq l}$. Hence {\rm
(\ref{jetid})} appears as a special case of {\rm (\ref{eqparam})},
without an additional formal parameter
$t$.}
\end{Rem}
For any ideal $I\subset \C\dbl y\dbr$, $y\in \C^r$, and any
 nonnegative integers $k,l$, with $k\geq 1$, we define an
ideal
$I^{(l)}\subset \C
\dbl
\Lambda_0\dbr [\hat \Lambda]$, where $\Lambda=(\Lambda_{0},
\hat \Lambda)$ are coordinates on $J^l_0(\C^k,\C^r)$ as in
(\ref{coord}) and (\ref{Lambda}), as follows:
\begin{multline}\Label{IL}
I^{(l)}:= \big\{ h\in \C \dbl \Lambda_0\dbr
[\hat \Lambda]: h(j_x^lF)=0\ {\rm for\ all}\ F:(\C^k_x,0)\to
(\C^r_y,0)\\ {\rm such\  that}\  u\circ F=0,\ {\rm for\
all}\ u
\in I\big\}.
\end{multline}
We have the following proposition.

\begin{Pro}\Label{awful}
If $I\subset \C \dbl y\dbr$ is a manifold ideal of
codimension $d$, then the ideal 
$I^{(l)}\subset \C
\dbl
\Lambda_0\dbr [\hat \Lambda]$ defined by {\rm (\ref{IL})}
is also a manifold ideal. Moreover, if the manifold ideal
$I$ is generated by
$\rho_1(y),\ldots,\rho_d(y)$ in $\C \dbl y\dbr$, then
the ideal $I^{(l)}$ in $\C
\dbl
\Lambda_0\dbr [\hat \Lambda]$ is generated by the components of
$\rho_1^{(l)}(\Lambda),\ldots,\rho_d^{(l)}(\Lambda)$ , where $\rho_j^{(l)}$ is
given by Proposition {\rm \ref{gener}}.
\end{Pro}
\begin{proof}
Recall by Proposition \ref{gener} that an invertible formal
map
$\psi :(\C^r_y,0)\to (\C^r_{y'},0)$ induces a formal
invertible map $\psi^{(l)}: (J_0^l(\C^k,\C^r),0)\to
(J_0^l(\C^k,\C^r),0)$. We leave it to the reader to check that
the equality 
\begin{equation}\Label{eqalide} 
(\psi^{(l)})_{*}(I^{(l)})=(\psi_{*}(I))^{(l)}
\end{equation}
follows from Proposition \ref{gener}, where the pushforward
of an ideal is given by (\ref{push}). If
$\rho_1(y),\ldots,\rho_d(y)$ are generators of the manifold
ideal $I$, we may choose a formal invertible map
$\psi:(\C^r_y,0)\to (\C^r_{y'},0)$ such  that
$y_j'=\psi_j(y)=\rho_j(y)$ for $j=1,\ldots,d$ and hence the
manifold ideal $\psi_{*}(I)\subset \C \dbl y'\dbr$ is
generated by the coordinate functions $y_1',\ldots,y_d'$.
We take $\Lambda=(\Lambda_0,\hat \Lambda)$ for coordinates
in the source jet space $J_0^l(\C^k,\C^r)$ and
$\Lambda'=(\Lambda'_0,\hat
\Lambda')$ for coordinates in the target one, as in (\ref{coord}) and
(\ref{Lambda}). It is then easy to check that the ideal
$(\psi_{*}(I))^{(l)}\subset
\C
\dbl \Lambda'_0 \dbr [\hat \Lambda']$ is the manifold ideal
generated by the coordinate functions
$(\Lambda'_{\alpha,i})$ for $0\leq |\alpha|\leq l$ and
$i=1,\ldots,d$. It follows from (\ref{eqalide}) that
$I^{(l)}\subset \C\dbl \Lambda_0\dbr [\hat \Lambda]$ is a
manifold ideal and is generated by the
$\psi^{(l)}_{\alpha,i}(\Lambda)$ for $0\leq |\alpha|\leq
l$ and
$i=1,\ldots,d$. Since by construction
$\psi_i(y)=\rho_i(y)$, $i=1,\ldots,d$, the last
part of the proposition follows.
\end{proof}

\def\rotrh{\hbox{\hskip1pt $\hskip-2.5pt\tilde\rho'{}\hskip-1pt^H$}}
\def\rorh{\hbox{\hskip1pt $\hskip-2.5pt\rho'{}\hskip-1pt^H$}}
\def\rotlh{\hbox{\hskip1pt $^H\hskip-2.5pt\tilde\rho'$}}
\def\rolh{\hbox{\hskip1pt $^H\hskip-2.5pt\rho'$}}
\def\rolho{\hbox{\hskip1pt $^{H^0}\hskip-3.5pt\rho'$}}
\def\rolhc{\hbox{\hskip1pt $^{\check
H}\hskip-3pt\check \rho'$}}

\section{Generators of the ideal ${\mathcal
I}(\scrM')^{(l)}$}\Label{IM} 
In this section, we consider a formal
generic manifold
${\scrM'}\subset \C^{N'}_{Z'}\times \C^{N'}_{\zeta'}$ of
codimension $d'$. Let
$$\rho':(\C^{N'}_{Z'}\times
\C^{N'}_{\zeta'},0)\to
(\C^{d'},0)$$
be a formal mapping
such that ${\mathcal
I}(\scrM')=(\rho'
(Z',\zeta'))=(\rho'_1(Z',\zeta'),\ldots,\rho'_{d'}
(Z',\zeta'))$ in $\C \dbl Z',\zeta'\dbr$. We define
\begin{equation}\Label{rho'}
{\tilde\rho}'(Z',\zeta'):= \overline{\rho'}(\zeta',Z').
\end{equation}
Since ${\mathcal I}(\scrM')$ is real,
the ideal ${\mathcal I}(\scrM')\subset \C\dbl Z',\zeta'\dbr$
is also generated by the components of ${\tilde
\rho}'(Z',\zeta')$.

 Given a
formal map
$H:(\C^N_Z,0)\to (\C_{Z'}^{N'},0)$, we define
 two formal mappings $\rolh:(\C^{N}_{Z}\times
\C^{N'}_{\zeta'},0)\to
(\C^{d'},0)$ and $\rorh:(\C^{N'}_{Z'}\times
\C^N_{\zeta},0)\to (\C^{d'},0)$ as follows
\begin{equation}\Label{reflec}
\rolh(Z,\zeta'):=\rho' (H(Z),\zeta'),\quad
\rorh(Z',\zeta):=\rho'
(Z',\overline{H}(\zeta)).
\end{equation}
Similarly, we define
\begin{equation}\Label{tildreflec}
\rotlh(Z,\zeta'):=\tilde\rho' (H(Z),\zeta'),\quad
\rotrh(Z',\zeta):=\tilde\rho'
(Z',\overline{H}(\zeta)).
\end{equation}
Note that by the reality condition, we have
\begin{equation}\Label{conjug} 
\rotlh(Z,\zeta')=\overline{\rorh}(\zeta',Z).
\end{equation}
Observe also that the components of
$\rolh(Z,\zeta')$ generate the reflection ideal
${\mathcal I}^H\subset
\C \dbl Z,\zeta'\dbr$ as defined by (\ref{crap}).

Throughout the rest of this section and \S\S
\ref{noM}--\ref{trick}, we fix
a nonnegative integer
$l$. Since $$\rho',\ \tilde \rho':(\C^{N'}_{Z'}\times
\C^{N'}_{\zeta'},0)\to (\C^{d'},0)$$
are formal mappings, by Proposition \ref{gener} there exist
unique formal mappings
\begin{equation}\Label{rhotilde}
\rho'{}^{(l)},\ {\tilde\rho}'{}^{(l)}:(J^l_0(\C^N,\C^{N'}
\times
\C^{N'}),0)\to (J^l_0(\C^N,\C^{d'}),0)
\end{equation}
 such that for every formal
mapping \begin{equation}\Label{F}
F=(F^1,F^2):(\C^N_{Z},0)\to
(\C^{N'}_{Z'}\times
\C^{N'}_{\zeta'},0),
\end{equation}
one has \begin{equation}\Label{psijet}
j_Z^l(\rho'\circ F)=\rho'{}^{(l)}(j_Z^lF),\quad 
j_Z^l(\tilde\rho'\circ F)={\tilde\rho}'{}^{(l)}(j_Z^lF).
\end{equation}
If $\Lambda=(\Lambda_{\alpha})_{|\alpha|\leq l}$, $\alpha
\in \N^N$, are the coordinates given
by (\ref{coord}) on the jet space
\begin{equation}\Label{split}
J^l_0(\C_Z^N,\C_{Z'}^{N'}\times
\C_{\zeta'}^{N'})=J^l_0(\C_Z^N,\C_{Z'}^{N'})\times 
J^l_0(\C_Z^N,\C_{\zeta'}^{N'}),
\end{equation} 
then we write 
$\Lambda=(\Lambda^1,\Lambda^2)$ according to the
splitting (\ref{split}). Thus, we have
$\Lambda^i=(\Lambda^i_{\alpha})_{|\alpha|\leq l}$, $i=1,2$, $\alpha
\in \N^N$.  As in (\ref{Lambda}), we continue to use the splitting
$\Lambda^i=(\Lambda_0^i,\hat \Lambda^i)$ with 
$\hat \Lambda^i=(\Lambda^i_{\alpha})_{1\leq |\alpha|\leq
l}$, $i=1,2$. 
Since $ I(\mathcal M')$ is generated either by the components
of
$\rho'(Z',\zeta')$ or by those of
$\tilde\rho'(Z',\zeta')$, it follows from Proposition
\ref{awful} that the ideal $ I(\mathcal M')^{(l)}$ in 
$\C\dbl\Lambda_0\dbr [\hat \Lambda]$ is
generated either by the components of 
$\rho'{}^{(l)}(\Lambda)$ or by the components of
$\tilde\rho'{}^{(l)}(\Lambda)$.

 We shall now give a more explicit expression for
$\rho'{}^{(l)}(\Lambda)$.  As in (\ref{crappy}), we write
$\rho'{}^{(l)}=(\rho_{\nu}'{}^{(l)})_{|\nu|\leq l}$ and
${\tilde\rho}'{}^{(l)}=({\tilde\rho}_{\nu}'{}^{(l)})
_{|\nu|\leq l}$. For any formal mapping $F(Z)$ as in
(\ref{F}), by (\ref{psijet}), the chain rule and
(\ref{reflec}), one has for any $\nu \in \N^N$,
$|\nu|\leq l$,

\begin{equation}\Label{pab}
\begin{aligned}
\rho_{\nu}'{}^{(l)}(j^l_ZF^1,j^l_ZF^2)&=
\frac{\partial^{\nu}}{\partial Z^{\nu}}\big[
\rho'(F^1(Z),F^2(Z))\big]=\frac{\partial^{\nu}}{\partial
Z^{\nu}}\big[\rho'{}^{\overline{F^2}}
(F^1(Z),Z)\big]\\
&= \sum_{\alpha \in \Nnn^{N'},\ \beta \in \Nnn^{N}\atop
|\beta|+|\alpha|\leq |\nu|,\, \beta \leq \nu}P_{\nu \alpha
\beta}\left((\partial^{\mu}F^1(Z))_{1\leq
|\mu|\leq
|\nu|}\right)\rho'{}^{\overline{F^2}}_{Z'{}^{\alpha}\zeta^{\beta}}
(F^1(Z),Z),
\end{aligned}
\end{equation}
where the $P_{\nu \alpha \beta}$ are universal scalar
polynomials depending only on $N$ and $N'$ (independent
of $F$ and
$\rho'$). Note that we also have 
\begin{equation}\Label{1}
P_{\nu0\nu}\equiv 1.
\end{equation}
 As in (\ref{reflec}), one should
regard
$\rho'{}^{\overline{F^2}}$ as a power series mapping of the
indeterminates
$(Z',\zeta)$; this is the meaning of the derivative
$\rho'{}^{\overline{F^2}}_{Z'{}^{\alpha}\zeta^{\beta}}$ in
(\ref{pab}). For any $\alpha \in \N^{N'}$, any
$\beta \in \N^{N}$ and for any formal map $F^2:(\C^N,0)\to
(\C^{N'},0)$, we have, again by the chain rule,
\begin{equation}\Label{rab}
\rho'{}^{\overline{F^2}}_{Z'{}^{\alpha}\zeta^{\beta}}(Z',\zeta)
=\sum_{\mu \in \Nn^{N'},\ |\mu|\leq |\beta|}R_{
\beta
\mu}\left((\partial^{\delta}F^2(\zeta))_{1\leq |\delta|\leq
|\beta|}\right)
\ \rho'_{Z'{}^{\alpha}\zeta'{}^{\mu}}(Z',F^2(\zeta)),
\end{equation}
where the  $R_{\beta \mu}$ are universal scalar polynomials
depending only on $N$ and $N'$ (independent of $F^2$ and
$\rho'$). Again, as in (\ref{rho'}), one should regard
$\rho'$ as a power series mapping of the
indeterminates $(Z',\zeta')$. Moreover, one has
$R_{00}=1$ and $R_{\beta0}=0$ for all $\beta\not =0$. As a
consequence of (\ref{pab}) and (\ref{rab}) and using the
notation (\ref{number}), we have for any formal mapping
$F(Z)$ as in (\ref{F}) 
\begin{multline}\Label{yt}
\rho_{\nu}'{}^{(l)}(j^l_ZF^1,j^l_ZF^2)=\\ \sum_{\alpha \in
\Nnn^{N'},\ \beta \in \Nnn^{N}\atop |\beta|+|\alpha|\leq
|\nu|,\, \beta \leq \nu}P_{\nu \alpha
\beta}\left(\hat \jmath_Z^lF^1\right)
\sum_{\mu \in \Nnn^{N'}\atop |\mu|\leq |\beta|}R_{
\beta
\mu}\left(\hat \jmath_Z^lF^2\right)
\ \rho'_{Z'{}^{\alpha}\zeta'{}^{\mu}}(F^1(Z),F^2(Z)).
\end{multline}
Hence, by the uniqueness in Proposition \ref{gener}, we
have in $\C \dbl \Lambda_0\dbr [\hat \Lambda]$,
$\Lambda=(\Lambda^1,\Lambda^2)$, for $\nu \in \N^N$,
$|\nu|\leq l$,
\begin{equation}\Label{rholambda1}
\rho'_{\nu}{}^{(l)}(\Lambda^1,\Lambda^2)=\sum_{\alpha \in
\Nnn^{N'},\ \beta \in \Nnn^{N}\atop |\beta|+|\alpha|\leq
|\nu|,\, \beta \leq \nu}P_{\nu \alpha
\beta}(\hat \Lambda^1)
\sum_{\mu \in \Nnn^{N'}\atop |\mu|\leq |\beta|}R_{
\beta
\mu}(\hat \Lambda^2)
\
\rho'_{Z'{}^{\alpha}\zeta'{}^{\mu}}(\Lambda_0^1,\Lambda_0^2).
\end{equation}
Using $^{F^1}\!\rho'(Z,F^2(Z))$ (given by (\ref{reflec}))
instead of
$\rho'{}^{\overline{F^2}} (F^1(Z),Z)$ in carrying out the
calculation in (\ref{pab}), one is led to the following
expression of $\rho'_{\nu}{}^{(l)}(\Lambda^1,\Lambda^2)$:

\begin{equation}\Label{rholambda2}
\rho'_{\nu}{}^{(l)}(\Lambda^1,\Lambda^2)=\sum_{\alpha \in
\Nnn^{N'},\ \beta \in \Nnn^{N}\atop |\beta|+|\alpha|\leq
|\nu|,\, \beta \leq \nu}P_{\nu \alpha
\beta}(\hat \Lambda^2)
\sum_{\mu \in \Nnn^{N'}\atop |\mu|\leq |\beta|}R_{
\beta
\mu}(\hat \Lambda^1)
\
\rho'_{Z'{}^{\mu}\zeta'{}^{\alpha}}
(\Lambda_0^1,\Lambda_0^2),
\end{equation}
where the polynomials $P_{\nu \alpha \beta}$ and $R_{\beta
\mu}$ are the same as those in (\ref{rholambda1}). Of
course, the expressions  (\ref{rholambda1}) and
(\ref{rholambda2}) also hold for $\rho'$ replaced by
$\tilde \rho'$ as well, since the components of $\tilde
\rho'$ are also generators of ${\mathcal I}(\scrM')$. 

We summarize the above in the following lemma.

\begin{Lem}\Label{cheney}
Let ${\scrM'}$, $\rho'$ and $\tilde\rho'$ be as above.
Then the ideal in 
$\C\dbl\Lambda_0\dbr [\hat \Lambda]$
generated by the components of $\rho'{}^{(l)}(\Lambda)$ is the
same as the ideal generated by the components of
${\tilde\rho}'{}^{(l)}(\Lambda)$, and both coincide with 
${\mathcal I}(\scrM')^{(l)}$.  Furthermore, the components
$\rho'_{\nu}{}^{(l)}(\Lambda)$ are given either
by $(\ref{rholambda1})$ or by $(\ref{rholambda2})$.
\end{Lem}

We
should mention that in what follows, we will use the
expression (\ref{rholambda1}) for
$\tilde\rho'$ and the expression (\ref{rholambda2}) for
$\rho'$, for a specific choice of $\rho'$.

\begin{Rem}\Label{RAB}
{\rm As in {\rm (\ref{rab})}, for any $\alpha \in
\N^{N'}$, any
$\beta \in \N^{N}$ and for any formal map $F^1:(\C^N,0)\to
(\C^{N'},0)$, we have
\begin{equation}\Label{rab2}
^{F^1}\!\!\rho'{}_{Z{}^{\beta}\zeta'{}^{\alpha}}(Z,\zeta')
=\sum_{\mu \in \Nn^{N'},\ |\mu|\leq |\beta|}R_{
\beta
\mu}\left((\partial^{\delta}F^1(Z))_{1\leq |\delta|\leq
|\beta|}\right)
\ \rho'_{Z'{}^{\mu}\zeta'{}^{\alpha}}(F^1(Z),\zeta'),
\end{equation}
where the universal polynomials $R_{\beta \mu}$ are the
same as those in {\rm (\ref{rab})}. Observe that
(\ref{rab2}) has already been used in proving
(\ref{rholambda2}).}
\end{Rem}

\section{Properties of reflection ideals and their
generators}\Label{noM}
As in \S \ref{IM}, we consider a formal generic
manifold 
$\scrM'\subset \C^{N'}_{Z'}\times \C^{N'}_{\zeta'}$ of
codimension $d'$. Since
${\scrM'}$ is generic, we may assume by using the formal implicit
function Theorem that 
$Z'=(z',w')\in
\C^{n'}\times
\C^{d'}$, $\zeta'=(\chi',\tau')\in \C^{n'}\times \C^{d'}$
with
$n'=N'-d'$, and that the ideal ${\mathcal
I}(\scrM')$ in $\C \dbl Z',\zeta'\dbr$ is given by
\begin{equation}\Label{Q'}
{\mathcal
I}(\scrM')=(w'-Q'(z',\zeta')),
\end{equation} 
where $Q':(\C^{n'}_{z'}\times \C^{N'}_{\zeta'},0)\to
(\C^{d'},0)$ is a formal mapping. Note that since 
${\scrM'}$ is real, we also have
\begin{equation}\Label{barQ'}
{\mathcal
I}(\scrM')=(\tau'-\bar{Q}'(\chi',Z')).
\end{equation}
For the rest of the paper, we make the following choice
of generators for ${\mathcal I}(\scrM')$
\begin{equation}\Label{bush}
\rho'(Z',\zeta'):=\tau'-\bar{Q}'(\chi',Z').
\end{equation}
Hence in view of (\ref{rho'}), we have
\begin{equation}\Label{gore}
\tilde\rho'(Z',\zeta')=w'-Q'(z',\zeta').
\end{equation}
We have the following proposition which holds for this
choice of generators of ${\mathcal I}(\scrM')$.
\begin{Pro}\Label{both}
Let ${\scrM'}\subset \C^{N'}_{Z'}
\times \C^{N'}_{\zeta'}$ be a formal generic manifold of
codimension $d'$ and
$H,H^0:(\C^N,0)\to (\C^{N'},0)$ be two formal mappings. Let
${\mathcal I}^H$ be the reflection ideal defined by {\rm
(\ref{crap})} and
$^H\!{\rho}'$ be the formal map given by {\rm
(\ref{reflec})} with the choice of $\rho'$ given by
{\rm (\ref{bush})}. Then the following hold.
\begin{enumerate}
\item [(i)] ${\mathcal I}^H={\mathcal I}^{H^0}$ $\iff$
$\rolh(Z,\zeta')=\rolho(Z,\zeta')$.
\item [(ii)] The reflection ideal ${\mathcal I}^H$ is
convergent (as in Definition {\rm \ref{pear}}) if and only
if  the components of
$\rolh(Z,\zeta')$ are convergent power series.
\item [(iii)] The reflection ideal ${\mathcal I}^H$ is
algebraic (as in Definition {\rm \ref{pear}}) if and only
if  the components of
$\rolh(Z,\zeta')$ are algebraic functions.
\end{enumerate}
\end{Pro}
\begin{proof}
(i) Since ${\mathcal I}^H=\left(\rolh(Z,\zeta')\right)$,
it follows that
$\rolh(Z,\zeta')=\rolho(Z,\zeta')$
implies the equality of the ideals 
${\mathcal I}^H$ and ${\mathcal I}^{H^0}$. Conversely, if
${\mathcal I}^H={\mathcal I}^{H^0}$, then there exists a
$d'\times d'$ matrix $a(Z,\zeta')$ with entries in $\C \dbl
Z,\zeta'\dbr$ such that 
\begin{equation}\Label{bof}
\rolh(Z,\zeta')=a(Z,\zeta')\, \rolho(Z,\zeta').
\end{equation}
Putting $\tau'=\bar{Q}'(\chi',H^{0}(Z))$ in (\ref{bof}) and
making use of (\ref{reflec}) and (\ref{bush}), we obtain
that
$\bar{Q}'(\chi',H^{0}(Z))=\bar{Q}'(\chi',H(Z))$ and hence 
$\rolh(Z,\zeta')=\rolho(Z,\zeta')$.

(ii) Since ${\mathcal I}^H=
\left(\rolh(Z,\zeta')\right)$, if the components of
$\rolh(Z,\zeta')$ are convergent, then
${\mathcal I}^H$ is convergent. Conversely, if ${\mathcal
I}^H$ is convergent, then by Definition \ref{pear} and Lemma
\ref{algebra} (ii), there exist $r_j(Z,\zeta')\in
\C\{Z,\zeta'\}$,
$j=1,\ldots,d'$, with linearly independent differentials
at 0
such that ${\mathcal I}^H=(r)=(r_1,\ldots,r_{d'})$ in $\C
\dbl Z,\zeta'\dbr$. As a consequence, there exist a
$d'\times d'$ invertible matrix $a(Z,\zeta')$ with entries
in 
$\C \dbl Z,\zeta'\dbr$ such that
\begin{equation}\Label{boffe}
r(Z,\zeta')=a(Z,\zeta')\rolh(Z,\zeta')=
a(Z,\zeta')\, (\tau'-\bar{Q}'(\chi',H(Z))),
\end{equation}
and hence $\partial
r/\partial \tau'(0)$ is invertible. By the implicit
function theorem, one sees that the equation
$r(Z,\chi',\tau')=0$ has a unique convergent solution
$\tau'=u(Z,\chi')$. It follows from (\ref{boffe}) that
$\bar{Q}'(\chi',H(Z))=u(Z,\chi')$ and hence that
$\rolh(Z,\zeta')$ is a convergent power series
mapping. This completes the proof of (ii).

(iii) The proof of this case is similar to that of part
(ii) above by making use of the algebraic version of
the implicit function theorem.
\end{proof}

\section{Ideals associated to formal generic manifolds and
mappings}\Label{new}
In this section, we consider two formal generic manifolds
$\scrM \subset \C^N_{Z}\times \C^{N}_{\zeta}$ and $\scrM'
\subset \C^{N'}_{Z'}\times \C^{N'}_{\zeta'}$ of
codimension $d$ and $d'$ respectively. 
We write $N=n+d$ and $N'=n'+d'$. As in \S \ref{noM},
we continue to use the choice of generators of ${\mathcal
I}(\scrM')$ given by (\ref{bush}) and (\ref{gore}). 

 By the implicit function theorem, there
exists a formal mapping 
\begin{equation}\Label{gamma}
\gamma:(\C^N_{\zeta}\times \C_t^n,0)\to (\C^N_Z,0),\quad
{\rm rk}\frac{\partial \gamma}{\partial t}(0)=n,
\end{equation}
such that for any $h\in {\mathcal I}(\scrM)$, 
$$h(\gamma
(\zeta,t),\zeta)=0,$$ and hence, by the reality of
${\scrM}$, we also have
\begin{equation}\Label{florida}
h(Z,\bar{\gamma}(Z,t))=0.
\end{equation}
Observe that each of the formal mappings $(\C^N\times
\C^n,0)
\ni (Z,t)\mapsto (Z,\bar{\gamma}(Z,t))$ and 
$(\C^N\times \C^n,0) \ni
(\zeta,t)\mapsto (\gamma(\zeta,t),\zeta)$ is a
parametrization of the formal generic manifold $\scrM$. If,
moreover,
${\scrM}\subset
\C^N\times
\C^N$ is the complexification of a generic real-analytic
(resp.\  real-algebraic) submanifold through the origin in
$\C^N$, then one can choose $\gamma$ to be convergent
(resp.\ algebraic). As in \cite{JAG}, we shall call a formal
map
$\gamma$ satisfying the above properties a {\em Segre
variety mapping} relative to
$\scrM$. Note that the formal map $(\C^n,0)\ni t\mapsto
\gamma (0,t)$ is a parametrization of ${\mathcal
S}_0(\scrM)$, the formal Segre variety of
$\scrM$ at 0 as defined in \S \ref{ideal}. In the rest of
this paper, we shall fix such a map $\gamma$. 

For a formal map $H:(\C^N,0)\to (\C^{N'},0)$ and the fixed
nonnegative integer $l$,
we define two formal mappings
$$\varphi^{[l]}(H;\cdot),\
\tilde\varphi^{[l]}(H;\cdot):\big(J^l_0(\C^N,
\C^{N'})\times \C^N_Z\times \C^n_{t},0\big)\to
J^l_0(\C^N, \C^{d'}),$$
as follows. Consider $\rho'{}^{(l)}(\Lambda^1,\Lambda^2)$ and
$\tilde\rho'{}^{(l)}(\Lambda^1,\Lambda^2)$ as defined in
(\ref{rhotilde}), with the choice of $\rho'$ and
$\tilde \rho'$ made in (\ref{bush}) and (\ref{gore}).  Taking
$\Lambda^2=\partial_Z^{\alpha}(\bar{H}(\bar{\gamma}(Z,t)))
_{|\alpha|\leq l}$ we set
\begin{equation}\Label{oregon}
\varphi^{[l]}(H;\Lambda^1,Z,t):=\rho'{}^{(l)}
\left(\Lambda^1,
(\partial_Z^{\alpha}(\bar{H}(\bar{\gamma}(Z,t)))
_{|\alpha|\leq l}\right),
\end{equation}
 and
\begin{equation}\Label{oregont}
\tilde\varphi^{[l]}(H;\Lambda^1,Z,t):=\tilde\rho'{}^{(l)}
\left(\Lambda^1,
(\partial_Z^{\alpha}(\bar{H}(\bar{\gamma}(Z,t)))
_{|\alpha|\leq l}\right).
\end{equation}
Observe that each component of the right hand side of
(\ref{oregon}) and (\ref{oregont}) is a formal power series
which is in
$\C\dbl\Lambda_0^1,Z,t\dbr [\hat
\Lambda^1]$. Here we recall that
$\Lambda^1=(\Lambda^1_0,\hat
\Lambda^1)$ are coordinates on $J^l_0(\C^N, \C^{N'})$ as in
(\ref{coord}) and (\ref{Lambda}).

 We shall
write, as in (\ref{crappy}),
$\varphi^{[l]}(H;\cdot)=(\varphi_{\nu}^{[l]}(H;\cdot))_
{|\nu|\leq l}$ and use a similar notation for 
$\tilde\varphi^{[l]}(H;\cdot)$. We shall now compute the
$\nu$-th component  $\tilde\varphi_{\nu}^{[l]}(H;\cdot)$.
It follows from (\ref{oregont}) and (\ref{rholambda1}), with
$\rho'$ replaced by $\tilde \rho'$, that

\begin{multline}\Label{seven}
\tilde\varphi_{\nu}^{[l]}(H;\Lambda^1,Z,t)=\\
\sum_{\alpha \in
\Nnn^{N'},\ \beta \in \Nnn^{N}\atop |\beta|+|\alpha|\leq
|\nu|,\, \beta \leq \nu}P_{\nu \alpha
\beta}(\hat \Lambda^1)
\sum_{\mu \in \Nnn^{N'}\atop |\mu|\leq |\beta|}R_{
\beta
\mu}(\hat \jmath_Z^l(\bar{H}(\bar{\gamma}(Z,t))))
\
\tilde\rho'_{Z'{}^{\alpha}\zeta'{}^{\mu}}(\Lambda_0^1,
\bar{H}(\bar{\gamma}(Z,t))),
\end{multline}
where $\hat
\jmath_Z^l(\bar{H}(\bar{\gamma}(Z,t)))=(\partial^{\delta}_Z
[\bar{H}(\bar{\gamma}(Z,t))])_{1\leq |\delta|\leq l}$.
By the chain rule, a computation similar to
(\ref{rab}) shows that one has
\begin{equation}\Label{calcul}
\frac{\partial^{|\beta|}}{\partial Z^{\beta}}\left[
\rotrh_{Z'{}^{\alpha}}(Z',\bar{\gamma}(Z,t))
\right]= \sum_{\mu \in \Nnn^{N'}\atop |\mu|\leq |\beta|}R_{
\beta
\mu}(\hat \jmath_Z^l(\bar{H}(\bar{\gamma}(Z,t))))
\
\tilde\rho'_{Z'{}^{\alpha}\zeta'{}^{\mu}}(Z',
\bar{H}(\bar{\gamma}(Z,t))),
\end{equation}
where the universal polynomials $R_{\beta \mu}$ are the
same as those in (\ref{rab}). On the other hand, by the
chain rule (again considering $\tilde\rho'{}^{H}$
as a power series mapping of the indeterminates
$(Z',\zeta)$), we also have 
\begin{equation}\Label{cac}
\frac{\partial^{|\beta|}}{\partial Z^{\beta}}\left[
\tilde\rho'{}^{\hskip-1pt H}_{Z'{}^{\alpha}}(Z',\bar{\gamma}(Z,t))
\right]=\sum_{|\delta|\leq |\beta|}c_{\beta \delta}(Z,t)
\tilde\rho'{}^{\hskip-1pt H}_{Z'{}^{\alpha}\zeta{}^{\delta}}
(Z',\bar{\gamma} (Z,t)).
\end{equation}
Here, the formal power series maps $c_{\beta
\delta}:(\C^N_Z\times \C^n_{t},0)\to \C$ depend only on
the Segre variety mapping $\gamma$ and not on the mapping
$H$. Moreover, if
$\gamma$ is convergent (resp.\  algebraic), then the
$c_{\beta \delta}$ are also convergent (resp.\  algebraic).
As a consequence of (\ref{seven}), (\ref{calcul}) and
(\ref{cac}), we obtain 
\begin{equation}\Label{systtilde}
\tilde\varphi_{\nu}^{[l]}(H;\Lambda^1,Z,t)=
\sum_{\alpha \in
\Nnn^{N'},\ \beta \in \Nnn^{N}\atop |\beta|+|\alpha|\leq
|\nu|,\, \beta \leq \nu}P_{\nu \alpha
\beta}(\hat \Lambda^1)\sum_{|\delta|\leq |\beta|}c_{\beta \delta}(Z,t)
\tilde\rho'{}^{\hskip-1pt H}_{Z'{}^{\alpha}\zeta{}^{\delta}}
(\Lambda_0^1,\bar{\gamma}
(Z,t)).
\end{equation}

If $H:(\C_Z^N,0)\to (\C_{Z'}^{N'},0)$ is a formal map such
that its complexification ${\mathcal H}:(\C^N_Z\times
\C^N_{\zeta},0)\to (\C^{N'}_{Z'}\times \C^{N'}_{\zeta'},0)$
given by (\ref{H}) sends $\scrM$ into ${\scrM'}$, then it
follows from (\ref{florida}) that
\begin{equation}\Label{maps}
h'(H(Z),\bar{H}(\bar{\gamma} (Z,t))) =0,\quad \forall
h'\in {\mathcal I}({\scrM'}). 
\end{equation}
Taking $h'$ in (\ref{maps}) to be any of the
components of $\rho'(Z',\zeta')$ or 
$\tilde\rho'(Z',\zeta')$ and making use of 
Remark
\ref{parameter}, we obtain
\begin{equation}\Label{lieb}
\rho'{}^{(l)}(j_Z^lH(Z),j^l_{Z}(\bar{H}(\bar{\gamma}(Z,
t))))=0,\quad
\tilde\rho'{}^{(l)}(j_Z^lH(Z),j^l_{Z}(\bar{H}(\bar{\gamma}(Z,
t))))=0.
\end{equation} 
Hence $\Lambda^1=j_Z^lH$ is a solution of each of
the systems of equations 
\begin{equation}\Label{ou}
\varphi^{[l]}(H;\Lambda^1,Z,t)=0,\quad
\tilde\varphi^{[l]}(H;\Lambda^1,Z,t)=0,
\end{equation}
where $\varphi^{[l]}(H;\cdot)$ and
$\tilde\varphi^{[l]}(H;\cdot)$ are defined by
(\ref{oregon}) and (\ref{oregont}) respectively.

We summarize the above in the following lemma.

\begin{Lem}\Label{summary} Let $H:(\C_Z^N,0)\to
(\C_{Z'}^{N'},0)$ and $\varphi^{[l]}(H;\cdot)$ and
$\tilde\varphi^{[l]}(H;\cdot)$ be the formal series given by  
$(\ref{oregon})$ and $(\ref{oregont})$ respectively.  Then the
ideal in
$\C\dbl\Lambda_0^1,Z,t\dbr [\hat
\Lambda^1]$ generated  by the
components of $\varphi^{[l]}(H;\Lambda^1,Z,t)$ is the
same as that generated by the components of
$\tilde\varphi^{[l]}(H;\Lambda^1,Z,t)$.  Moreover,
the components of $\tilde\varphi^{[l]}(H;\Lambda^1,Z,t)$ are
given by formula $(\ref{systtilde})$.  If, in addition, the
complexification $\mathcal H$ of $H$, as given by $(\ref
{H})$, maps
$\mathcal M$ into
$\mathcal M'$, then $\Lambda^1=j_Z^lH$ is a solution of each of
the two systems of equations in $(\ref{ou})$.
\end{Lem}
\section{Iterated Segre mappings and associated
ideals}\Label{bushwin} In this section, we assume that
$\scrM$ and $\scrM'$ are given formal generic manifolds as
in \S \ref{new}. We continue to use the choice of
generators $\rho'(Z',\zeta')$ and $\tilde
\rho'(Z',\zeta')$ given in (\ref{bush}) and (\ref{gore})
for the ideal ${\mathcal I}(\scrM')$. If
$\gamma$ is a Segre variety mapping relative to $\scrM$ as
defined in (\ref{gamma}), we define, as in \cite{JAG}, the
{\em iterated Segre mappings} (relative to
$\scrM$) as follows. First, we set $v^0:=0\in \C^N$. For any
positive integer
$j$, $v^j:(\C^{nj},0)\to (\C^N,0)$ is the formal mapping
defined inductively by
\begin{equation}\Label{vj}
v^j(t^1,\ldots,t^j):=\gamma(\bar{v}^{j-1}(t^1,\ldots,t^{j-1}),
t^j),\quad t^1,\ldots,t^j\in \C^n.
\end{equation}
In what follows, it will be convenient to introduce for a
given positive integer $j$ the notation
$$t^{[j]}:=(t^1,\ldots,t^{j}),$$ considered as a variable in
$\C^{nj}$. With this notation, we may rewrite (\ref{vj})
in the form
$$v^j(t^{[j]})=\gamma(\bar{v}^{j-1}(t^{[j-1]}),t^j).$$
It follows from (\ref{florida}) and (\ref{vj}) that for any
$h\in {\mathcal I}(\scrM)$ and any nonnegative integer $j$,
we have
\begin{equation}\Label{fair}
h(v^j(t^{[j]}),\bar{v}^{j+1}(t^{[j+1]}))=0.
\end{equation}

If $H:(\C^N,0)\to (\C^{N'},0)$ is a
formal mapping and
$j$  is a fixed nonnegative integer,  we
define two formal mappings
$$\psi^{[l,j]}(H;\cdot),\ \tilde\psi^{[l,j]}(H;\cdot):
(J^l_0(\C^N,
\C^{N'})\times \C^{(j+2)n},0)\to
J^l_0(\C^N, \C^{d'}),$$
as follows:
\begin{equation}\Label{kansas}
\psi^{[l,j]}(H;\Lambda^1,t^{[j+2]}):=
\varphi^{[l]}(H;\Lambda^1,v^{j+1}(t^{[j+1]}),t^{j+2})
\end{equation}
and similarly
\begin{equation}\Label{alab}
\tilde\psi^{[l,j]}(H;\Lambda^1,t^{[j+2]}):=
\tilde\varphi^{[l]}(H;\Lambda^1,v^{j+1}
(t^{[j+1]}),t^{j+2}).
\end{equation}
Here we recall that the formal mappings
$\varphi^{[l]}(H;\cdot)$ and $\tilde\varphi^{[l]}(H;\cdot)$
are given by (\ref{oregon}) and  (\ref{oregont})
respectively. Hence the components of
$\psi^{[l,j]}(H;\cdot)$ and $\tilde\psi^{[l,j]}(H;\cdot)$
are formal power series in the ring $\C\dbl
\Lambda^1_0,t^{[j+2]}\dbr [\hat \Lambda^1]=\C
\dbl
\Lambda^1_0,t^1,\ldots,t^{j+2}\dbr[\hat\Lambda^1]$. It
follows from the definition of
$\tilde
\psi^{[l,j]}(H;\cdot)$ and from (\ref{systtilde}) that one
has the following identity for every $\nu\in \N^N$,
$|\nu|\leq l$,
\begin{multline}\Label{back}
\tilde
\psi_{\nu}^{[l,j]}(H;\Lambda^1,t^{[j+2]})=\\
\sum_{\alpha
\in
\Nnn^{N'},\ \beta \in \Nnn^{N}\atop |\beta|+|\alpha|\leq
|\nu|,\, \beta \leq \nu}P_{\nu \alpha
\beta}(\hat \Lambda^1)
\sum_{|\delta|\leq |\beta|}u^j_{\beta
\delta}(t^{[j+2]})
\tilde\rho'{}^{\hskip-1pt H}7_{Z'{}^{\alpha}\zeta{}^{\delta}}
(\Lambda_0^1,\bar{v}^{j+2}(t^{[j+2]})).
\end{multline}
Here we have used (\ref{vj}) in the form
$\bar{v}^{j+2}(t^{[j+2]})=\bar{\gamma}
(v^{j+1}(t^{[j+1]}),t^{j+2})$ and set $u^j_{\beta
\delta}(t^{[j+2]}):=c_{\beta
\delta}(v^{j+1}(t^{[j+1]}),t^{j+2})$, where the $c_{\beta
\delta}$ are as in (\ref{systtilde}). Note that since the
$c_{\beta \delta}$ are independent of $H$, so are the
$u_{\beta \delta}^j$.  Moreover, if
$\scrM$ is the complexification of a real-analytic (resp.\ 
real-algebraic) generic submanifold of
$\C^N$ through the origin, then we may assume that the
formal power series 
$u^j_{\beta \delta}$ in (\ref{back}) are convergent (resp.\ 
algebraic). The
following lemma is then a consequence of Lemmas
\ref{cheney} and
\ref{summary} as well as the above
construction.

\begin{Lem}\Label{add}
The ideals
$\left(\psi^{[l,j]}(H;\Lambda^1,t^{[j+2]})\right)$ and
$\left(\tilde
\psi^{[l,j]}(H;\Lambda^1,t^{[j+2]})\right)$ in $\C\dbl
\Lambda^1_0,t^{[j+2]}\dbr [\hat \Lambda^1]$ are the
same. In particular, let
$$S:(\C^{n(j+2)},0)\to J^l_0(\C^N,\C^{N'}),$$
be a formal map with
$S(t^{[j+2]})=(S_0(t^{[j+2]}),\hat S(t^{[j+2]}))$ as in
{\rm (\ref{Lambda})} and
$S_0(0)=0$. Then
$\Lambda^1=S(t^{[j+2]})$ is a solution of  
$\psi^{[l,j]}(H;\Lambda^1,t^{[j+2]})=0$ if
and only if it is a solution of
$\tilde\psi^{[l,j]}(H;\Lambda^1,t^{[j+2]})=0$.
Moreover, if
$H:(\C^N,0)\to (\C^{N'},0)$ is a formal map such
that its complexification ${\mathcal H}$
sends $\scrM$ into ${\scrM'}$, then
\begin{equation}\Label{put}
\Lambda^1=
\left((\partial^{\alpha}H)(v^{j+1}(t^{[j+1]}))\right)_{|\alpha|\leq
l}
\end{equation}
is a solution of the systems of equations 
\begin{equation}\Label{systeq}
\psi^{[l,j]}(H;\Lambda^1,t^{[j+2]})=0,\quad \tilde
\psi^{[l,j]}(H;\Lambda^1,t^{[j+2]})=0.
\end{equation}
\end{Lem}

We need the following lemma concerning the iterated Segre
mappings of ${\scrM}$.
\begin{Lem}\Label{segre}
Let $\gamma (\zeta,t)$ be a Segre variety mapping relative
to the generic formal manifold ${\scrM}\subset \C^N\times
\C^N$ as defined in {\rm (\ref{gamma})} and $v^j$ the
iterated Segre mappings as defined in {\rm (\ref{vj})}.
Then for every nonnegative integer $j$, there exists a
unique formal mapping $\xi^{j}:(\C^{n(j+1)},0)\to
(\C^n,0)$ such that
\begin{equation}\Label{switch}
v^{j+2}(t^{[j+1]},\xi^j(t^{[j+1]}))=v^j(t^{[j]}).
\end{equation}
Moreover, if the formal mapping $\gamma$ is convergent
(resp.\ algebraic), then $\xi^j$ is convergent (resp.\
algebraic).
\end{Lem}
\begin{proof}
Since $\scrM$ is generic, by making use of the implicit
function theorem we can assume that
$Z=(z,w)\in
\C^n\times \C^d$, where $d$ is the codimension of $\scrM$
and $n=N-d$, and that ${\mathcal I}(\scrM)$ is generated by
the components of $w-Q(z,\zeta)$ where $Q:(\C^{n+N},0)\to
(\C^d,0)$ is a formal mapping. If $\scrM$ is the
complexification of a real-analytic (resp.\ real-algebraic)
generic submanifold $M\subset \C^N$, then the formal map $Q$
is convergent (resp.\ algebraic). By the reality of $\scrM$,
one has
\begin{equation}\Label{real}
Q(z,\chi,\bar{Q}(\chi,z,w))=w.
\end{equation}
Corresponding to the splitting $Z=(z,w)$, we may write
$\gamma(\zeta,t)=(\mu(\zeta,t),\nu(\zeta,t))$, with
$\mu:(\C^{N+n},0)\to (\C^n,0)$ and $\nu :(\C^{N+n},0)\to
(\C^d,0)$. By the definition of a Segre variety mapping
$\gamma$, we necessarily have 
\begin{equation}\Label{tri}
\gamma(\zeta,t)=(\mu(\zeta,t),Q(\mu (\zeta,t),\zeta)).
\end{equation}
Since $\rk \, \partial \gamma/\partial t(0)=n$,  the matrix
$\partial \mu/\partial t(0)$ is invertible. As a
consequence of the implicit function theorem, there exist a
formal mapping $\pi :(\C^{N+n},0)\to (\C^n,0)$ such that 
\begin{equation}\Label{drink}
\mu (\bar{\gamma}(Z,t),\pi
(Z,t))=z,\quad {\rm where}\ 
Z=(z,w).
\end{equation} 
It follows from (\ref{real}), (\ref{tri}) and
(\ref{drink}) that
$\gamma (\bar{\gamma}(Z,t),\pi
(Z,t))=Z$. The lemma follows by taking
$\xi^j(t^{[j+1]}):=\pi
(v^{j}(t^{[j]}),t^{j+1})$.\begin{footnote} {If one takes $\mu
(\zeta,t)=t$, where $\mu
(\zeta,t)$ is the component of $\gamma(\zeta,t)$ as in
(\ref{tri}), then the reader can check that one has
$\xi^j(t^{[j+1]})=t^j$.   }
\end{footnote}
\end{proof}

If $j$ is a nonnegative integer, $l$ is the previously
fixed nonnegative integer, and
$H:(\C^N,0)\to (\C^{N'},0)$ is a formal map, then for
$\nu\in \N^N$,
$|\nu|\leq l$,  and
$\epsilon\in \N^n$, we define formal mappings
\begin{equation}\Label{bad}
\Theta^{[l,j]}_{\nu,\epsilon}(H;\cdot),\ 
\tilde\Theta^{[l,j]}_{\nu,\epsilon}(H;\cdot):
(J^l_0(\C^N,\C^{N'})\times
\C^{n(j+1)},0)\to \C^{d'}
\end{equation}
by
\begin{equation}\Label{theta}
\begin{aligned}
\Theta^{[l,j]}_{\nu,\epsilon}(H;\Lambda^1,t^{[j+1]}):&=
\partial^{\epsilon}_{t^{j+2}}\psi^{[l,j]}_{\nu}
(H;\Lambda^1,t^{[j+2]})\big|_{t^{j+2}=\xi^{j}(t^{[j+1]})},\\
\tilde\Theta^{[l,j]}_{\nu,\epsilon}(H;\Lambda^1,t^{[j+1]}):&=
\partial^{\epsilon}_{t^{j+2}}\tilde\psi^{[l,j]}_{\nu}
(H;\Lambda^1,t^{[j+2]})\big|_{t^{j+2}=\xi^{j}(t^{[j+1]})},
\end{aligned}
\end{equation}
where $\psi^{[l,j]}(H;\cdot)=(\psi^{[l,j]}
_{\nu}(H;\cdot))_{|\nu|\leq l}$ and 
$\tilde\psi^{[l,j]}(H;\cdot)=(\tilde\psi^{[l,j]}
_{\nu}(H;\cdot))_{|\nu|\leq l}$
are defined by (\ref{kansas}) and (\ref{alab}) respectively,
and the map
$\xi^j$ is given by Lemma
\ref{segre}. Observe that each component of
$\Theta^{[l,j]}_{\nu,\epsilon}(H;\cdot)$ and 
$\tilde\Theta^{[l,j]}_{\nu,\epsilon}(H;\cdot)$ is a formal
power series in $\C \dbl
\Lambda^1_0,t^{[j+1]}\dbr[\hat\Lambda^1]$. We have the
following lemma concerning the formal power series
mapping $\tilde\Theta^{[l,j]}_{\nu,\epsilon}(H;\cdot)$.

\begin{Lem}\Label{last}
For any $\nu\in \N^N$, $|\nu|\leq l$ and any $\epsilon \in
\N^n$, the following holds.
\begin{equation}\Label{least}
\tilde\Theta^{[l,j]}_{\nu,\epsilon}(H;\Lambda^1,t^{[j+1]})=
\sum_{|\alpha|\leq l \atop |\delta|\leq
l+|\epsilon|}\omega^j_{\nu \epsilon \alpha
\delta}(\hat\Lambda^1,t^{[j+1]})\tilde\rho'{}^{H}_{Z'{}^
{\alpha}\zeta^{\delta}}(\Lambda^1_0,\bar{v}^j(t^{[j]})),
\end{equation}
where each $\omega^j_{\nu \epsilon \alpha
\delta}(\hat \Lambda^1,t^{[j+1]})\in \C\dbl
t^{[j+1]}\dbr[\hat\Lambda^1]$ is independent of the formal
mapping $H$. Here, $\tilde\rho'{}^{H}(Z',\zeta)$ is
considered as a formal power series mapping in the indeterminates
$(Z',\zeta)$. Moreover, if the Segre variety mapping
$\gamma$ relative to
$\scrM$ is convergent  $($resp.\  algebraic$)$, then each
formal power series $\omega^j_{\nu \epsilon \alpha
\delta}(\hat \Lambda^1,t^{[j+1]})$ is in $\C\{
t^{[j+1]}\}[\hat\Lambda^1]$ {\rm (}resp.\ ${\mathcal A}\{
t^{[j+1]}\}[\hat\Lambda^1]${\rm )}.
\end{Lem}
\begin{proof}
The proof is an immediate consequence of (\ref{back}), the
definition of the 
$\tilde\Theta^{[l,j]}_{\nu,\epsilon}(H;\cdot)$ given in
(\ref{theta}), and (\ref{switch}).
\end{proof} 

The formal power series given by (\ref{theta}) will not be
used until
\S 10.  Their importance lies in the following remark.
\begin{Rem}\Label{S}
{\rm Let $S:(\C^{n(j+1)},0)\to J_0^l(\C^N,\C^{N'})$ be a
formal mapping such that
$S_0(0)=0$ where
$S(t^{[j+1]})=(S_{\nu}(t^{[j+1]}))_{|\nu|\leq
l}=(S_0(t^{[j+1]}),\hat S(t^{[j+1]}))$ as in (\ref{Lambda}).
Then,
$\Lambda^1=S(t^{[j+1]})$ is a solution of the system of
equations $\tilde\psi^{[l,j]}(H;\Lambda^1,t^{[j+2]})=0$ if
and only if it is a solution of the system of equations 
$\tilde\Theta^{[l,j]}_{\nu,\epsilon}(H;\Lambda^1,t^{[j+1]})
=0$ for all $\nu \in \N^N$, $|\nu|\le l$, and all
$\epsilon \in \N^n$. This is an immediate consequence
of the fact that $S(t^{[j+1]})$ is independent of the
indeterminate $t^{j+2}$ and the definition (\ref{theta}) of
the $\tilde\Theta^{[l,j]}_{\nu,\epsilon}(H;\cdot)$.}
\end{Rem}
\section{Properties of solutions of the system $\psi
^{[l,j]}(H;\Lambda^1,t^{[j+2]})=0$}\Label{trick}
The following technical lemma will be essential for the
proofs of Theorems \ref{fdi}, \ref{convideal} and
\ref{algideal}.
\begin{Lem}\Label{court}
Let $\scrM\subset \C^N\times \C^N$ and
$\scrM'\subset \C^{N'}\times \C^{N'}$ be
formal generic manifolds and $H:(\C^N,0)\to
(\C^{N'},0)$ a formal map such that its
complexification ${\mathcal H}:(\C^N\times
\C^N,0)\to (\C^{N'}\times \C^{N'},0)$
 sends $\scrM$ into ${\scrM'}$.
Assume that $H$ is not totally degenerate as in
Definition {\rm \ref{apple}}. Let $l,j$ be nonnegative
integers and
$\psi^{[l,j]}(H;\cdot)$ the formal map given by {\rm
(\ref{kansas})}. Let $S:(\C^{n(j+1)},0)\to J_0^l(\C^N,\C^{N'})$,
$S_0(0)=0$, be a formal map and assume that
$\Lambda^1=S(t^{[j+1]})=(S_{\nu}(t^{[j+1]}))_{|\nu|\leq
l}=(S_0(t^{[j+1]}),\hat S(t^{[j+1]}))$ is a formal solution
of the system
\begin{equation}\Label{systpsi}
\psi^{[l,j]}(H;\Lambda^1,t^{[j+2]})=0.
\end{equation}
Then the following holds. For every $\nu \in \N^N$,
$|\nu|\leq l$,
\begin{equation}\Label{key}
\rolh(v^{j+1}(t^{[j+1]}),\zeta')=
\sum_{|\mu|\leq |\nu|}R_{\nu\mu}(\hat
S(t^{[j+1]}))\rho'_{Z'{}^{\mu}}(S_0(t^{[j+1]}),\zeta'),
\end{equation}
where $\rho'$ and $\rolh$ are given by
{\rm (\ref{bush})} and {\rm (\ref{reflec})} respectively,
and the $R_{\nu\mu}$ are the universal polynomials given in
{\rm (\ref{rab})}. Here
$\rho'(Z',\zeta')$ and $\rolh(Z,\zeta')$ are
considered as formal power series mappings in the
indeterminates
$(Z',\zeta')$ and $(Z,\zeta')$ respectively. If, moreover,
$S(t^{[j+1]})=\left((\partial^{\alpha}H^0)
(v^{j+1}(t^{[j+1]}))\right)_{|\alpha|\leq l}$ for some
formal map
$H^0:(\C^N,0)\to (\C^{N'},0)$, then {\rm (\ref{key})} for
$|\nu|\leq l$ is equivalent to
\begin{equation}\Label{flue}
\rolh_{Z^{\nu}}(v^{j+1}(t^{[j+1]}),\zeta')=\rolho_{Z^{\nu}}
(v^{j+1}(t^{[j+1]}),\zeta').
\end{equation}
\end{Lem}
\begin{proof}
In what follows, we use the coordinates $Z'=(z',w')$,
$\zeta'=(\chi',\tau')$ as in the beginning of \S \ref{noM}
and write
\begin{equation}\Label{decomp}
H(Z)=(f(Z),g(Z)),\quad  {\rm
with}
\ z'=f(Z)\ {\rm and} \ w'=g(Z).
\end{equation} 
For the proof of (\ref{key}), we proceed by induction on
$|\nu|$ and we start first by proving (\ref{key})
for
$\nu =0$. Note that since
$\Lambda^1=S(t^{[j+1]})$ is a solution  of the system
(\ref{systpsi}), it follows that 
$\psi_{0}^{[l,j]}(H;S(t^{[j+1]}),t^{[j+2]})=0$. The latter
equation is equivalent to
\begin{equation}\Label{refer1}
\rho'(S_0(t^{[j+1]}),\bar{H}(\bar{v}^{j+2}(t^{[j+2]})))=0.
\end{equation}
Observe that since ${\mathcal H}$ maps $\scrM$ into
$\scrM'$, we have by making use of (\ref{maps}) and
(\ref{vj}) that 
\begin{equation}\Label{refer2}
\rho'(H(v^{j+1}(t^{[j+1]})),\bar{H}(\bar{v}^{j+2}(t^{[j+2]})))=0.
\end{equation}
It follows from (\ref{decomp}), (\ref{refer1}),
(\ref{refer2}) and (\ref{bush}) that 
\begin{equation}\Label{plug}
\bar{Q}'\left(\bar{f}(\bar{v}^{j+2}(t^{[j+2]})),
S_0(t^{[j+1]})\right)=
\bar{Q}'\left(\bar{f}(\bar{v}^{j+2}(t^{[j+2]})),
H(v^{j+1}(t^{[j+1]}))\right).
\end{equation}
To show that (\ref{key}) holds for $\nu=0$, in view of
(\ref{reflec}) and (\ref{bush}), we must show  that 
\begin{equation}\Label{display}
\bar{Q}'\left(\chi',
S_0(t^{[j+1]})\right)=
\bar{Q}'\left(\chi',
H(v^{j+1}(t^{[j+1]}))\right).
\end{equation}
For this, by e.g.\ Proposition 5.3.5 of \cite{BER}, it
suffices to show that $\Rk B$, the rank of the formal map
\begin{equation}\Label{B} B:(\C^{n(j+2)}_{t^{[j+2]}},0)\to
(\C^{n(j+1)}\times
\C^{n'},0)
\end{equation} given by
$B(t^{[j+2]}):=(t^{[j+1]},\bar{f}(\bar{v}^{[j+2]}
(t^{[j+2]})))$, is $n(j+1)+n'$. The latter 
follows from the fact that $H$ is not totally degenerate.
Indeed, since $v^1(t)=\gamma (0,t)$ is a parametrization
of the formal Segre variety ${S}_0(\scrM)$, it follows from
Definition
\ref{apple} that
$\Rk H\circ v^1=n'$ . Hence we also have $\Rk f\circ
v^1=n'$. From this, we easily obtain that $\Rk
B=n(j+1)+n'$. This completes the proof of (\ref{key}) for
$\nu=0$.

It follows from (\ref{rholambda2}) and the definition of the
$\psi^{[l,j]}(H;\cdot)$ given in (\ref{kansas}) that the
following identity holds for all $\nu \in \N^N$,
$|\nu|\leq l$,
\begin{multline}\Label{formula1}
\psi^{[l,j]}_{\nu}(H;\Lambda^1,t^{[j+2]})=\\
\sum_{\alpha \in
\Nnn^{N'},\ \beta \in \Nnn^{N}\atop |\beta|+|\alpha|\leq
|\nu|,\, \beta \leq \nu}P_{\nu \alpha
\beta}(\widehat m(t^{[j+2]}))
\sum_{\mu \in \Nnn^{N'}\atop |\mu|\leq |\beta|}R_{
\beta
\mu}(\hat \Lambda^1)
\
\rho'_{Z'{}^{\mu}\zeta'{}^{\alpha}}
(\Lambda_0^1,m_0(t^{[j+2]})),
\end{multline}
where
\begin{equation}\Label{m}
\begin{aligned}
m_{\alpha}(t^{[j+2]}):&=\partial_Z^{\alpha}\left[\bar{H}
(\bar{\gamma}(Z,t))\right]\Big|_{Z=v^{j+1}(t^{[j+1]}),
\, t=t^{j+2}}\ ,\\
\widehat m
(t^{[j+2]}):&=\left(m_{\alpha}(t^{[j+2]})\right)_{1\leq
|\alpha|\leq l}\ .
\end{aligned}
\end{equation}
In view of (\ref{1}), we may rewrite
(\ref{formula1}) as follows
\begin{multline}\Label{formula2}
\psi^{[l,j]}_{\nu}(H;\Lambda^1,t^{[j+2]})=
\sum_{\mu \in \Nnn^{N'}\atop |\mu|\leq |\nu|}R_{
\nu
\mu}(\hat \Lambda^1)
\
\rho'_{Z'{}^{\mu}}
(\Lambda_0^1,m_0(t^{[j+2]}))+\\
\sum_{\alpha \in
\Nnn^{N'},\ \beta \in \Nnn^{N}\atop |\beta|+|\alpha|\leq
|\nu|,\, \beta<\nu}P_{\nu \alpha
\beta}(\widehat m (t^{[j+2]}))
\sum_{\mu \in \Nnn^{N'}\atop |\mu|\leq |\beta|}R_{
\beta
\mu}(\hat \Lambda^1)
\
\rho'_{Z'{}^{\mu}\zeta'{}^{\alpha}}
(\Lambda_0^1,m_0(t^{[j+2]}))
\end{multline}

 Let $\epsilon \in \N^N$, $0<|\epsilon|\leq l$, and assume
that (\ref{key}) holds for all $\nu \in \N^N$ with
$|\nu|<|\epsilon|$. We now show that (\ref{key}) holds
for $\nu=\epsilon$. Since
$\Lambda^1=S(t^{[j+1]})$ is a  solution of
the system (\ref{systpsi}), it follows from
(\ref{formula2}), with $\nu$ replaced by $\epsilon$, that we
have
\begin{multline}\Label{formula3}
\sum_{\mu \in \Nnn^{N'}\atop |\mu|\leq |\epsilon|}R_{
\epsilon
\mu}(\hat S(t^{[j+1]}))
\
\rho'_{Z'{}^{\mu}}
(S_0(t^{[j+1]}),m_0(t^{[j+2]})))=\\
-\sum_{ |\beta|+|\alpha|\leq
|\epsilon|\atop \beta<\epsilon}P_{\epsilon \alpha
\beta}( \widehat m(t^{[j+2]}))
\sum_{\mu \in \Nnn^{N'}\atop |\mu|\leq |\beta|}R_{
\beta
\mu}(\hat S(t^{[j+1]}))
\
\rho'_{Z'{}^{\mu}\zeta'{}^{\alpha}}
(S_0(t^{[j+1]}),m_0(t^{[j+2]}))).
\end{multline}
On the other hand, using the notation 
\begin{equation}\Label{ee}
e_{\alpha}(t^{[j+1]}):=(\partial^{\alpha}H)(v^{j+1}
(t^{[j+1]})),\quad
\hat e(t^{[j+1]}):=(e_{\alpha}(t^{[j+1]}))_{1\leq
|\alpha|\leq l},
\end{equation}
it follows from Lemma \ref{add} that
$\Lambda^1=(e_{\alpha}(t^{[j+1]}))_{|\alpha|\leq l}$ is
also a solution of (\ref{systpsi}) and hence, from
(\ref{formula2}), we obtain
\begin{multline}\Label{formula4}
\sum_{\mu \in \Nnn^{N'}\atop |\mu|\leq |\epsilon|}R_{
\epsilon
\mu}(\hat e (t^{[j+1]}))
\
\rho'_{Z'{}^{\mu}}
(e_0(t^{[j+1]}),m_0(t^{[j+2]})))=\\
-\sum_{ |\beta|+|\alpha|\leq
|\epsilon|\atop \beta<\epsilon}P_{\epsilon \alpha
\beta}( \widehat m (t^{[j+2]}))
\sum_{\mu \in \Nnn^{N'}\atop |\mu|\leq |\beta|}R_{
\beta
\mu}(\hat e(t^{[j+1]}))
\
\rho'_{Z'{}^{\mu}\zeta'{}^{\alpha}}
(e_0(t^{[j+1]}),m_0(t^{[j+2]}))).
\end{multline}
By (\ref{rab2}) with $F^1=H$, $Z$ replaced by
$v^{j+1}(t^{[j+1]})$ and $\zeta'$ replaced by
$m_0(t^{[j+2]})$, we have for any $\beta\in \N^N$,
$|\beta|\leq l$ and any $\alpha \in \N^{N'}$,
\begin{multline}\Label{replace}
\rolh_{Z{}^{\beta}\zeta'{}^{\alpha}}
(v^{j+1}(t^{[j+1]}),
m_0(t^{[j+2]}))
=\\
\sum_{\mu \in \Nn^{N'},\ |\mu|\leq |\beta|}R_{
\beta
\mu}\left(\hat e(t^{[j+1]})\right)
\
\rho'_{Z'{}^{\mu}\zeta'{}^{\alpha}}
(e_0(t^{[j+1]}),m_0(t^{[j+2]})).
\end{multline}
By the induction hypothesis, since $\beta <\epsilon$ in the
right hand side of (\ref{formula4}), we have, after
differentiating (\ref{key}) (with
$\nu=\beta$) with respect to $\zeta'$ and replacing
$\zeta'$ by $m_0(t^{[j+2]})$,
\begin{multline}\Label{replace2}
\rolh_{Z{}^{\beta}\zeta'{}^{\alpha}}
(v^{j+1}(t^{[j+1]}),
m_0(t^{[j+2]}))
=\\
\sum_{\mu \in \Nn^{N'},\ |\mu|\leq |\beta|}R_{
\beta
\mu}\left(\hat S(t^{[j+1]})\right)
\
\rho'_{Z'{}^{\mu}\zeta'{}^{\alpha}}
(S_0(t^{[j+1]}),m_0(t^{[j+2]})).
\end{multline}
It follows from (\ref{formula3}), (\ref{formula4}),
(\ref{replace}) and (\ref{replace2}) that 

\begin{multline}\Label{fof}
\sum_{\mu \in \Nnn^{N'}\atop |\mu|\leq |\epsilon|}R_{
\epsilon
\mu}(\hat S(t^{[j+1]}))
\rho'_{Z'{}^{\mu}}
(S_0(t^{[j+1]}),m_0(t^{[j+2]})))=\\
\sum_{\mu \in
\Nnn^{N'}\atop |\mu|\leq |\epsilon|}R_{
\epsilon
\mu}(\hat e (t^{[j+1]}))
\
\rho'_{Z'{}^{\mu}}
(e_0(t^{[j+1]}),m_0(t^{[j+2]}))).
\end{multline}
Using (\ref{replace}) with $\beta =\epsilon$ and $\alpha
=0$, we obtain that (\ref{fof}) implies 
\begin{multline}\Label{foffe}
\rolh_{Z{}^{\epsilon}}
(v^{j+1}(t^{[j+1]}),
m_0(t^{[j+2]}))=
\sum_{\mu \in \Nnn^{N'}\atop |\mu|\leq |\epsilon|}R_{
\epsilon
\mu}(\hat S(t^{[j+1]}))
\rho'_{Z'{}^{\mu}}
(S_0(t^{[j+1]}),m_0(t^{[j+2]}))).
\end{multline}
To prove (\ref{key}) for $\nu =\epsilon$, we must show that
(\ref{foffe}) still holds if
$m_0(t^{[j+2]})=\bar{H}(\bar{v}^{j+2}(t^{[j+2]}))$ is
replaced by an arbitrary $\zeta'=(\chi',\tau')\in \C^{N'}$.
Observe that for $\mu
\in \N^{N'}$, $|\mu|>0$, we have in view of (\ref{bush}), 
$\rho'_{Z'{}^{\mu}}(Z',\zeta')=-
\bar{Q}_{Z'{}^{\mu}}'(\chi',Z'):=a_{\mu}(Z',\chi')$ and
since
$|\epsilon|>0$, 
$\rolh_{Z^{\epsilon}}(Z,\zeta')=-
\partial_Z^{\epsilon}\left[
\bar{Q}'(\chi',H(Z))\right]:=b_{\epsilon}(Z,\chi')$. Recall
also that since
$\epsilon \not =0$, $R_{\epsilon 0}=0$ (see \S \ref{IM}).
Hence, (\ref{foffe}) may be rewritten in the form
\begin{multline}\Label{rue}
b_{\epsilon}(v^{j+1}(t^{[j+1]}),\bar{f}(\bar{v}^{j+2}
(t^{[j+2]})))=\\
\sum_{\mu \in \Nnn^{N'}\atop 0<|\mu|\leq
|\epsilon|}R_{
\epsilon
\mu}(\hat S(t^{[j+1]})) a_{\mu}(S_0(t^{[j+1]}),
\bar{f}(\bar{v}^{j+2}
(t^{[j+2]}))),
\end{multline}
and we must show that (\ref{rue}) still holds with $\bar{f}(\bar{v}^{j+2}
(t^{[j+2]}))$ replaced by an arbitrary $\chi'\in
\C^{n'}$. For this, one can apply the same rank argument
using the map $B$ defined in (\ref{B}) as in the case
$\nu=0$. This completes the proof of ({\ref{key}}). 

To complete the proof of Lemma \ref{court}, it suffices to
observe that the equivalence of (\ref{key}) and
(\ref{flue}) follows from (\ref{rab2}). 
\end{proof}

\section{Proof of Theorem \ref{fdi}}\Label{close}
For the proof of Theorem \ref{fdi}, we shall need 
Proposition \ref{bla} given below. We assume that $\scrM$,
$\scrM'$ and the iterated Segre mappings $v^j$ are as in \S
\ref{bushwin} and continue to use the notation of that
section. In particular, we still assume that
$\rho'(Z',\zeta')$ and $\tilde\rho'(Z',\zeta')$ are the
special choice of generators of ${\mathcal I}(\scrM')$
given by (\ref{bush}) and (\ref{gore}). 
\begin{Pro}\Label{bla}
Let $\scrM \subset \C^N\times \C^N$ and $\scrM'\subset
\C^{N'}\times \C^{N'}$ be formal generic manifolds. Let
$H^0:(\C^N,0)\to (\C^{N'},0)$ be a formal mapping such its
complexification
${\mathcal H}^0$ sends $\scrM$ into $\scrM'$. Assume that
$H^0$ is not totally degenerate $($as in
Definition {\rm \ref{apple}}$)$. Then for every pair of
nonnegative integers
$l,j$, there exists a positive integer $K=K(H^0,l,j)$ such
that if
$H:(\C^N,0)\to (\C^{N'},0)$ is a formal map whose
complexification ${\mathcal H}$ maps $\scrM$ into $\scrM'$
and such that
\begin{equation}\Label{hypo}
\rolho_{Z^{\delta}}(v^j(t^{[j]}),\zeta')=\rolh_{Z^{\delta}}(v^j(t^{[j]}),\zeta'),\quad
|\delta|\leq K,
\end{equation}
then
\begin{equation}\Label{result}
\rolho_{Z^{\delta}}(v^{j+1}(t^{[j+1]}),\zeta')=\rolh_{Z^{\delta}}(v^{j+1}(t^{[j+1]}),\zeta'),\quad
|\delta|\leq l.
\end{equation}
Here $\rolh(Z,\zeta')$ and 
$\rolho(Z,\zeta')$ are the formal mappings given by
{\rm (\ref{reflec})} with the choice {\rm (\ref{bush})} of
$\rho'$.
\end{Pro}
\begin{proof} We fix the pair of nonnegative integers
$l,j$. In the ring
$R:=\C
\dbl
 \Lambda_0^1, t^{[j+1]}\dbr[\hat \Lambda^1]$, where
$\Lambda^1=(\Lambda^1_0,\hat \Lambda^1)$ are coordinates
on $J^l_0(\C^N,\C^{N'})$ as in (\ref{coord}) and (\ref{Lambda}), we
consider the ideal
${\mathcal J}$ generated by the components of the formal
mappings
$$\tilde\Theta^{[l,j]}_{\nu,\epsilon}(H^0;\Lambda^1,t^{[j+1]}),
\quad \nu \in \N^N,\ |\nu|\leq l,\ \epsilon \in \N^n,$$ 
where the
$\tilde\Theta^{[l,j]}_{\nu,\epsilon}(H^0;\cdot)$ are given
by (\ref{theta}). Since $R$ is Noetherian, there exists a
positive integer 
$L=L(H^0,l,j)$ such that the ideal ${\mathcal J}$ is
generated by the components of the formal mappings
$$\tilde\Theta^{[l,j]}_{\nu,\epsilon}(H^0;\Lambda^1,t^{[j+1]}),
\quad \nu \in \N^N,\ |\nu|\leq l,\ \epsilon
\in \N^n,\ |\epsilon|\leq L.$$ 

We claim that the conclusion of Proposition \ref{bla} holds
with $K:=L+l$. Indeed, let $H:(\C^N,0)\to (\C^{N'},0)$ be a
formal map whose complexification sends $\scrM$ into
$\scrM'$ and such that (\ref{hypo}) holds (with this choice
of $K$). We must prove that (\ref{result}) holds. By
(\ref{conjug}), we have
\begin{equation}\Label{plus}
\rotrh(Z',\zeta)=\overline{\rolh}(\zeta,Z'),\quad
\tilde\rho'{}^{\hskip-1pt H^0}(Z',\zeta)=\overline{\rolho}(\zeta,Z'),
\end{equation}
and hence it follows from (\ref{hypo}) that
\begin{equation}\Label{oui}
\tilde\rho'{}^{\hskip-1pt H^0}_{Z'{}^{\alpha}
\zeta^{\delta}}(Z',\bar{v}^j(t^{[j]}))=
\tilde\rho'{}^{\hskip-1pt H}_{Z'{}^{\alpha}\zeta^{\delta}}
(Z',\bar{v}^j(t^{[j]})),\quad \alpha \in \N^{N'},\
|\delta|\leq K,
\end{equation}
where we have considered $\rotrh(Z',\zeta)$ and
$\tilde\rho'{}^{\hskip-1pt H^0}(Z',\zeta)$ as formal mappings in the
indeterminates $(Z',\zeta)$ as in (\ref{tildreflec}). As a
consequence of (\ref{oui}), (\ref{least}) and the choice of
$K$, it follows that
\begin{equation}\Label{state}
\tilde\Theta^{[l,j]}_{\nu,\epsilon}
(H;\Lambda^1,t^{[j+1]})=\tilde\Theta^{[l,j]}
_{\nu,\epsilon}(H^0;\Lambda^1,t^{[j+1]}),\quad |\nu|\leq
l,\quad |\epsilon|\leq L.
\end{equation}
By Lemma \ref{add},
\begin{equation}\Label{ridi}
\Lambda^1=\left((\partial^{\alpha}H)(v^{j+1}(t^{[j+1]}))
\right)_{|\alpha|\leq l}
\end{equation} is a formal solution of the system of
equations $\tilde\psi^{[l,j]}(H;\Lambda^1,t^{[j+2]})=0$, and
hence by Remark
\ref{S} (since (\ref{ridi}) is independent of the
indeterminate $t^{j+2}$), it is also a solution of the
system of equations
$\tilde\Theta^{[l,j]}_{\nu,\epsilon}
(H;\Lambda^1,t^{[j+1]})=0$ for $|\nu|\leq l$ and all
$\epsilon \in \N^n$. From (\ref{state}), we conclude that
(\ref{ridi}) is also a solution of the system
of
equations
$$\tilde\Theta^{[l,j]}
_{\nu,\epsilon}(H^0;\Lambda^1,t^{[j+1]})=0,\quad |\nu|\leq
l,\quad |\epsilon|\leq L.$$
By the choice of $L$, it follows that the formal
power series mapping given by (\ref{ridi}) is a solution of
the system of equations
$$\tilde\Theta^{[l,j]}
_{\nu,\epsilon}(H^0;\Lambda^1,t^{[j+1]})=0,\quad |\nu|\leq
l,\quad \forall \epsilon \in \N^n.$$
Again making use of Remark \ref{S}, we conclude that
(\ref{ridi}) is a formal solution of the system of
equations $\tilde\psi^{[l,j]}(H^0;\Lambda^1,t^{[j+2]})=0$
and hence, by Lemma \ref{add}, also a solution of the system
of equations $\psi^{[l,j]}(H^0;\Lambda^1,t^{[j+2]})=0$. We
may now apply Lemma
\ref{court} with
$H$ and
$H^0$ interchanged and with
$S(t^{[j+1]})=\left((\partial^{\alpha}H)(v^{j+1}(t^{[j+1]}))
\right)_{|\alpha|\leq l}$ to conclude that (\ref{flue})
holds, which is the desired conclusion
(\ref{result}) of Proposition \ref{bla}.\end{proof}

\begin{proof}[Proof of Theorem {\rm \ref{fdi}}]
Since $\scrM$ is of finite type, it follows from Theorem 2.3
in
\cite{JAG} (see also \cite{MA}) and the definition of
finite type given in
\S
\ref{ideal}, that there exists an integer
$j_0$,
$2\le j_0\le d+1$, where $d$ is the codimension of $\scrM$
such that
$\Rk  v^{j_0}=N$. By applying Proposition \ref{bla} $j_0$
times, we conclude that there exists an integer
$K_0=K_0(H^0)>0$
\footnote{To find $K_0$, we proceed as follows. We define
inductively a finite sequence of nonnegative integers 
$K_{q}$, $0\le q\le j_0$, by putting $K_{j_0}=0$ and
$K_{q}=K(H^0,K_{q+1},q)$ where $K(H^0,l,j)$ is the integer
given by Proposition \ref{bla}. } such that if
$H:(\C^N,0)\to (\C^{N'},0)$ is a formal map whose
complexification
${\mathcal H}$ sends
$\scrM$ into $\scrM'$ and such that
\begin{equation}\Label{hypo2}
\rolho_{Z^{\delta}}(v^0,\zeta')=\rolh_{Z^{\delta}}(v^0,\zeta'),\quad
|\delta|\leq K_0,
\end{equation}
then
\begin{equation}\Label{hypo3}
\rolho(v^{j_0}(t^{[j_0]}),\zeta')=\ 
\rolh(v^{j_0}(t^{[j_0]}),\zeta').
\end{equation}
Recall that $v^0=0\in \C^N,$ and hence we may rewrite
(\ref{hypo2}) in the form
\begin{equation}\Label{partial}
\partial_Z^{\delta}\left[\rho'(H(Z),\zeta')\right]\big|_{Z=0}
=\partial_Z^{\delta}\left[\rho'(H^0(Z),\zeta')
\right]\big|_{Z=0},\quad |\delta|\leq K_0.
\end{equation}
It is then clear that if $H$ is
a formal map such that ${\mathcal H}$ sends $\scrM$ 
into $\scrM'$ with
$j_0^{K_0}H=j_{0}^{K_0}H^0$, then (\ref{partial}) and hence
(\ref{hypo2}) and (\ref{hypo3}) hold. Since $\Rk
v^{j_0}=N$, it follows e.g.\ from Proposition 5.3.5 of
\cite{BER} that (\ref{hypo3}) implies
\begin{equation}\Label{no}
\rolho(Z,\zeta')=\rolh(Z,\zeta').
\end{equation}
{}From the definition of the reflection ideal ${\mathcal
I}^H$ given in (\ref{crap}), we conclude that (\ref{no})
implies that the reflection ideals
${\mathcal I}^H$ and ${\mathcal I}^{H^0}$ are the same. The
proof of Theorem \ref{fdi} is complete.
\end{proof}

\section{Proof  of Theorem \ref{convideal}}
In this section, we consider two germs $(M,0)$ and
$(M',0)$ of
real-analytic generic submanifolds in $\C^N$ and $\C^{N'}$
respectively. We let $\scrM\subset \C^{N}_{Z}\times
\C^{N}_{\zeta}$ and
$\scrM'\subset \C^{N'}_{Z'}\times \C^{N'}_{\zeta'}$ be their
complexifications. For generators of ${\mathcal
I}(\scrM')$, we take a convergent mapping $\rho'(Z',\zeta')$
as in (\ref{bush}). We shall also use the
corresponding notation for $\tilde\rho'(Z',\zeta')$ given
by (\ref{gore}). Moreover, we choose a convergent Segre
variety mapping
$\gamma$ relative to ${\scrM}$ as defined in (\ref{gamma});
hence the corresponding iterated Segre mappings $v^{j}$
defined in (\ref{vj}) are also convergent. Using the
notation of \S \ref{bushwin}, we have the following
proposition.

\begin{Pro}\Label{final}
Let $(M,0)$ and $(M',0)$ be germs of generic real-analytic
submanifolds in $\C^N$ and $\C^{N'}$ of
codimension $d$ and $d'$ respectively. Let
$H:(\C^N,0)\to (\C^{N'},0)$ be a formal map sending $M$
into $M'$. Assume that $H$ is not totally degenerate $($as
 in Definition {\rm \ref{apple}}$)$. Then for every
nonnegative integer $j$, the following holds. If
\begin{equation}\Label{U}
\rolh_{Z^{\beta}}(v^j(t^{[j]}),\zeta')\in
(\C\{t^{[j]},\zeta'\})^{d'},\quad \forall \beta \in \N^{N},
\end{equation}
then 
\begin{equation}\Label{V}
\rolh_{Z^{\beta}}(v^{j+1}(t^{[j+1]}),\zeta')\in
(\C\{t^{[j+1]},\zeta'\})^{d'},\quad \forall \beta \in
\N^{N}.
\end{equation}
Here $^H\!\!\rho'(Z,\zeta')$ is the formal 
mapping given by {\rm (\ref{reflec})} relative to the
choice of the convergent mapping $\rho'(Z',\zeta')$ given
by {\rm (\ref{bush})}.
\end{Pro}
\begin{proof}
We fix a pair of nonnegative integers $l,j$, and we shall
prove that if (\ref{U}) holds,
then 
\begin{equation}\Label{WW}
\rolh_{Z^{\nu}}(v^{j+1}(t^{[j+1]}),\zeta')\in
(\C\{t^{[j+1]},\zeta'\})^{d'},\quad \forall \nu \in
\N^{N},\ |\nu|\le l.
\end{equation}
The proposition will clearly follow. 

It follows from (\ref{U}) and (\ref{conjug}) that one has
\begin{equation}\Label{UB}
\tilde\rho'{}^{\hskip-1pt H}_{\zeta^{\nu}}(Z',\bar{v}^j(t^{[j]}))\in
(\C\{Z',t^{[j]}\})^{d'},\quad \forall \nu \in \N^{N},
\end{equation}
where $\rotrh(Z',\zeta)$ is the formal mapping
given by (\ref{tildreflec}). It follows from Lemma
\ref{last} and (\ref{UB}) that the components of the formal
power series mappings
$\tilde\Theta^{[l,j]}_{\nu,\epsilon}(H;\Lambda^1,t^{[j+1]})$,
 for $\nu \in \N^N$, $|\nu|\leq l$,
and $\epsilon \in \N^n$, (defined by (\ref{theta})), are in
the ring $\C \{\Lambda^1_0,t^{[j+1]}\}[\hat \Lambda^1]$.
(We should observe at this point that the components of
the formal mappings
$\tilde\psi^{[l,j]}(H;\Lambda^1,t^{[j+2]})$, defined in
(\ref{alab}), are not yet known to be convergent.) By Lemma
\ref{add}, it follows that 
$\Lambda^1=\left((\partial^{\alpha}H)(v^{j+1}(t^{[j+1]}))
\right)_{|\alpha|\leq l}$ is a formal solution of the system
of equations
\begin{equation}\Label{meal}
\tilde\psi^{[l,j]}(H;\Lambda^1,t^{[j+2]})=0,
\end{equation}
and hence by Remark \ref{S}, it is also a formal solution
of the system of equations 
\begin{equation}\Label{artin}
\tilde\Theta^{[l,j]}_
{\nu,\epsilon}(H;\Lambda^1,t^{[j+1]})=0,\quad |\nu|\leq l,\
\epsilon \in \N^n.
\end{equation}
Since the mappings $\tilde\Theta^{[l,j]}_
{\nu,\epsilon}(H;\cdot)$ are convergent, it follows from
Artin's approximation theorem \cite{artin1} that there
exists a convergent solution of (\ref{artin}) given by
$\Lambda^1=S(t^{[j+1]})=(S_0(t^{[j+1]}),\hat
S(t^{[j+1]}))$, where
\begin{equation}\Label{salah}
S:(\C^{n(j+1)},0)\to
J^l_0(\C^N,\C^{N'}),\quad S(0)=j_0^lH.
\end{equation}
Since the convergent mapping
$S(t^{[j+1]})$ is independent of the variable $t^{j+2}$, it
follows from Remark \ref{S} that $\Lambda^1=S(t^{[j+1]})$
is also a solution of the system of equations given by
(\ref{meal}). Hence, by Lemma \ref{add},
$\Lambda^1=S(t^{[j+1]})$ is a solution of the system of
equations
$$\psi^{[l,j]}(H;\Lambda^1,t^{[j+2]})=0.$$
We may now apply Lemma \ref{court} for the convergent
solution $S(t^{[j+1]})$ to obtain (\ref{key}). To conclude
that (\ref{WW}) holds, it suffices to observe that the
right hand side of (\ref{key}) is a convergent map. This
completes the proof of Proposition \ref{final}.
\end{proof}

\begin{proof}[Proof of Theorem {\rm \ref{convideal}}]
Since $M$ is of finite type at 0, by Theorem 10.5.5 of
\cite{BER} (see also \cite{Acta,JAG}),  there exists
an integer
$k_0$,
$2\leq k_0\le 2(d+1)$ (where
$d$ is the codimension of $M$) such that in any neighborhood
$U$ of $0\in \C^{nk_0}$, there exists $t^{[k_0]}_0\in U$
such that 
\begin{equation}\Label{meet}
\rk \frac{\partial
v^{k_0}}{\partial t^{[k_0]}}(t^{[k_0]}_0)=N,\quad
v^{k_0}(t^{[k_0]}_0)=0.
\end{equation}
Since $v^0=0\in \C^N$, we observe that for any multiindex
$\beta\in
\N^N$,
\begin{equation}\Label{silen}
\rolh_{Z^{\beta}}(v^0,\zeta')=
\partial_{Z}^{\beta}\big[\rho'(H(Z),\zeta')\big]\big|_{Z=0}\in
(\C\{\zeta'\})^{d'},\quad \forall \beta \in \N^{N}.
\end{equation}
Applying Proposition \ref{final} $k_0$ times, we conclude
in particular that
\begin{equation}\Label{sky}
\rolh(v^{k_0}(t^{[k_0]}),\zeta')\in
(\C\{t^{[k_0]},\zeta'\})^{d'}.
\end{equation}
Hence there exists an open neighborhood $U\times
V\subset \C^{nk_0}\times \C^{N'}$ of 0 where the mapping 
$\rolh(v^{k_0}(t^{[k_0]}),\zeta')$ is convergent. If
we choose $t^{k_0}_0\in U$ such that (\ref{meet}) holds and
apply the rank theorem, we obtain that the mapping 
$\rolh(Z,\zeta')$ is convergent. By the definition of
the reflection ideal ${\mathcal I}^H$ given in (\ref{crap})
and Definition
\ref{pear}, it follows that ${\mathcal
I}^H$ is convergent. This completes the proof of Theorem
\ref{convideal}.
\end{proof}

\begin{Rem}
{\rm As mentioned in \S \ref{int}, Theorem \ref{convideal}
was first proved in \cite{Mirh} for an invertible formal
map
$H$ and in the case where $M$ and $M'$ are real-analytic
hypersurfaces in
$\C^N$. We should point out here that the techniques used
in this paper are somewhat different from those of
\cite{Mirh}. For instance, the use of Cauchy estimates 
was a crucial tool in \cite{Mirh}, but is not needed in our
approach in this paper. We should also note that Corollary
7.4 and Theorem 7.1 in \cite{Mirh}, which are proved there
in the case of invertible formal mappings between
real-analytic hypersurfaces of finite type, can be extended
to the case of finite formal mappings between generic
real-analytic submanifolds of finite type of
$\C^N$ by making use of Theorem \ref{convideal}. We do not
give any further details.}
\end{Rem}
\section{Proof of Theorem \ref{algideal}}\Label{2.7}

In this section, we consider two germs $(M,0)$ and
$(M',0)$ of
real-algebraic generic submanifolds in $\C^N$ and $\C^{N'}$
respectively. We let $\scrM\subset \C^{N}_{Z}\times
\C^{N}_{\zeta}$ and
$\scrM'\subset \C^{N'}_{Z'}\times \C^{N'}_{\zeta'}$ be their
complexifications. For generators of ${\mathcal
I}(\scrM')$, we take the components of an algebraic mapping
$\rho'(Z',\zeta')$ as in (\ref{bush}). We also  use the
corresponding notation for $\tilde\rho'(Z',\zeta')$ given
by (\ref{gore}). Moreover, we choose an algebraic Segre
variety
 mapping $\gamma$ relative to ${\scrM}$ as defined in
(\ref{gamma}); hence the corresponding iterated Segre
mappings $v^{j}$ defined in (\ref{vj}) are also algebraic.
We have the following analog of Theorem \ref{convideal} for
generic real-algebraic submanifolds.

\begin{Thm}\Label{copy}
Let $(M,0)$ and $(M',0)$ be germs of real-algebraic generic
submanifolds in $\C^N$ and $\C^{N'}$ respectively and
$H:(\C^N,0)\to (\C^{N'},0)$ a formal map sending $M$
into $M'$. Assume that $M$ is of finite type at $0$ and $H$
is not totally degenerate. Then the reflection ideal
${\mathcal I}^H$, as defined by {\rm (\ref{crap})}, is
algebraic.
\end{Thm}
This theorem will be used in the proof of Theorem
\ref{algideal} in the case where $H$ is a convergent
mapping. The proof of Theorem \ref{copy} follows the same
lines as that of Theorem \ref{convideal}, by making use of
the following analog of Proposition \ref{final} in the
algebraic setting. 

\begin{Pro}\Label{finalalg}
Let $(M,0)$ and $(M',0)$ be germs of generic real-algebraic
submanifolds in $\C^N$ and $\C^{N'}$ of
codimension $d$ and $d'$ respectively. Let
$H:(\C^N,0)\to (\C^{N'},0)$ be a formal map sending $M$
into $M'$. Assume that $H$ is not totally degenerate
$($as in Definition {\rm \ref{apple}}$)$. Then for
every nonnegative integer
$j$, the following holds. If
\begin{equation}\Label{if}
\rolh_{Z^{\beta}}(v^j(t^{[j]}),\zeta')\in
({\mathcal A}\{t^{[j]},\zeta'\})^{d'},\quad \forall \beta
\in
\N^{N},
\end{equation}
then
\begin{equation}\Label{then}
\rolh_{Z^{\beta}}(v^{j+1}(t^{[j+1]}),\zeta')\in
({\mathcal A}\{t^{[j+1]},\zeta'\})^{d'},\quad \forall \beta
\in
\N^{N}.
\end{equation}
Here $\rolh(Z,\zeta')$ is the formal 
mapping given by {\rm (\ref{reflec})} relative to the
choice of the algebraic mapping $\rho'(Z',\zeta')$ given
by {\rm (\ref{bush})}.
\end{Pro}
\begin{proof} 
The proof of this proposition follows very closely that of
Proposition \ref{final}. One has  to note that all the
convergent mappings involved in the latter are also
algebraic in the present case. Also, the convergent
solution $S(t^{[j+1]})$ of the system (\ref{artin}),
 given in (\ref{salah}) and obtained by
making use of Artin's approximation theorem, can be chosen
to be algebraic. Indeed, in the present case, the mappings
involved in (\ref{artin}) are algebraic and another version
of Artin's approximation theorem \cite{artin2} yields a
solution which is also algebraic. We omit further details.
\end{proof}

\begin{proof}[Proof of Theorem {\rm \ref{algideal}}]
Choose $U,U'\subset \C^N$ two open polydiscs centered at
 the origin such that $H$ is holomorphic in $U$ and
$H(U\cap M)\subset U'\cap M'$. We may assume that the
real-algebraic generic submanifold $M'\subset \C_{Z'}^N$ is
given by $\tilde \rho'(Z',\bar{Z}')=0$ where
\begin{equation}
\tilde \rho'(Z',\bar{Z}'):= w'-Q'(z',\bar{Z}'),\quad
Z'=(z',w')\in \C^n\times \C^d,
\end{equation}
with $\tilde \rho'(Z',\zeta')$ a $\C^d$ valued
algebraic map defined in $U'\times U'$. Here we
recall that $d$ is the codimension of $M$ (and of $M'$) and
$n=N-d$. Equivalently, $M'$ is also given by $
\rho'(Z',\bar{Z}')=0$ where
\begin{equation}\Label{GTD}
 \rho'(Z',\bar{Z}'):=
\bar{w}'-\bar{Q}'(\bar{z}',Z').
\end{equation}
To prove Theorem \ref{algideal}, by Proposition \ref{both}
(iii), it suffices to show that the convergent generators
$\rolh(Z,\zeta')$ of the reflection ideal ${\mathcal
I}^H$ are algebraic, where we have used the notation given
by (\ref{crap}) and (\ref{reflec}). Since the Jacobian of
$H$ is not identically zero and there is no germ at 0 of a
nonconstant  holomorphic function
$h:(\C^N,0)\to \C$ with
$h(M)\subset \bR$, it follows that there
exists
$p_0\in U\cap M$ such that $M$ is of finite type at $p_0$
and the Jacobian of $H$ at $p_0$ is not zero (see
e.g Lemma 13.3.2 of \cite{BER}). Put
$p'_0:=H(p_0)\in U'\cap M'$. We define the translation maps
$\varphi_{p_0}(Z):=Z-p_0$ and $\varphi_{p'_0}(Z'):=Z'-p'_0$.
We put $M_{p_0}:=\varphi_{p_0}(M)$ and
$M'_{p'_0}:=\varphi_{p'_0}(M')$. Observe that
$M_{p_0}$ and $M'_{p'_0}$ are real-algebraic generic
submanifolds through the origin in
$\C^{N}$ 
with $M_{p_0}$ of finite type at 0. We also define
\begin{equation}\Label{view}
\check H(\check
Z):=(\varphi_{p'_0}\circ H\circ
\varphi_{p_0}^{-1})(\check Z)
\end{equation} for $\check Z$ close enough
to the origin in $\C^{N}$. We can regard $\check H$ as a germ
at the origin of a biholomorphism  sending the
germ $(M_{p_0},0)$ onto $(M'_{p'_0},0)$. Note also that the
germ
$(M'_{p'_0},0)$ is defined by
$\check \rho'(\check Z',\overline{\check Z'})=0$ where
\begin{equation}\Label{hot}
\check \rho'(\check Z',{\check \zeta'}):=\check
\tau'+\bar{w}'_{p'_0}-\bar{Q}'(\check
\chi'+\bar{z}'_{p'_0},
\check Z'+p'_{0}),\quad \check \zeta'=(\check
\chi',\check \tau')\in \C^n\times \C^d,
\end{equation}
with $p'_0=(z'_{p'_0},w'_{p'_0})\in \C^n\times \C^d$. It
follows from Theorem \ref{copy} and Proposition \ref{both}
(iii) that the convergent mapping $\rolhc(\check Z,{\check
\zeta'})=\check \rho'(\check H(\check Z),\check \zeta')$ is
in $({\mathcal A}\{\check Z,\check \zeta\})^d$, i.e.\ that
the components of the map
$$(\C^N\times \C^n,0)\ni (\check Z,\check \chi')\mapsto
\bar{Q}'(\check
\chi'+\bar{z}'_{p'_0},
\check H(\check Z)+p'_{0})\in \C^d$$
are
in ${\mathcal A}\{\check Z,\check \zeta\}$. In view of
(\ref{view}), we conclude that the map
$$(\C^N\times \C^n,(p_0,\bar{z}'_{p'_0}))\ni (
Z,
\chi')\mapsto
\bar{Q}'(
\chi',
H(Z))\in \C^d$$
is algebraic i.e.\ each component of this map satisfies a
non-trivial polynomial equation with polynomial
coefficients for $Z$ near $p_0$ and $\chi'$ near
$\bar{z}'_{p'_0}$. By unique continuation, the same
equations hold for $(Z,\chi')$ close to $0\in  \C^{N+n}$.
This shows that the components of $\rolh(Z,\zeta')$ are
in
${\mathcal A}\{Z,\zeta'\}$ which gives the desired
conclusion of Theorem \ref{algideal}.
\end{proof}

\section{Proofs of Propositions \ref{link}, \ref{fdhol} and™
\ref{convhol} and Theorems \ref{plus1} and
\ref{plus2}}\Label{tenure}  In this section, we consider a
formal generic manifold
$\scrM'\subset \C^{N'}_{Z'}\times \C^{N'}_{\zeta'}$ of
codimension $d'$ and we assume that the ideal ${\mathcal
I}(\scrM')$ is generated by the components of the formal map
$\rho'(Z',\zeta')$ given by (\ref{bush}). We write
\begin{equation}\Label{game}
\rho'(Z',\zeta')=\tau'-\bar{Q}'(\chi',Z')=\tau'-
\sum_{\alpha \in \Nn^{n'}}q_{\alpha}(Z')\chi'{}
^{\alpha},
\end{equation}
where the
$q_{\alpha}(Z')=\big(q_{1,\alpha}(Z'),\ldots,
q_{d',\alpha}(Z')\big)$
are in
$(\C\dbl Z'\dbr)^{d'}$ and $n'=N'-d'$.

The proof of the
following criterion for holomorphic nondegeneracy of formal
generic manifolds is left to the reader (see e.g.\ \cite{St}
and \cite{BER}, Chapter 11,
for the case where
$\scrM'$ is the complexification of a real-analytic generic
submanifold).

\begin{Lem}\Label{gorelose}
The formal generic manifold $\scrM'$ as above is
holomorphically nondegenerate if and only if there exist
$\alpha^1,\ldots,\alpha^{N'}\in \N^{n'}$ and
$j_1,\ldots,j_{N'}\in \{1,\ldots,d'\}$ such that 
\begin{equation}\Label{det}
{\rm det}\left(\frac{\partial q_{j_l,\alpha^{l}}}{\partial
Z'_m}(Z')\right)_{1\le l,m\le N'}\not
= 0,\quad {\rm in}\ \C \dbl Z'\dbr,
\end{equation}
where the formal power series $q_{\alpha}(Z')$ are given by
{\rm (\ref{game})}.
\end{Lem}
We also need the following lemma for the proof of
Proposition \ref{fdhol}.
\begin{Lem}\Label{prop22} Let
$R(x,y)=(R_1(x,y),\ldots,R_r(x,y))\in (\C\dbl x,y\dbr)^r$,
$x\in
\C^q$, $y\in \C^r$, and $h^0:(\C^q,0)\to (\C^r,0)$ be a
formal map such that 
\begin{enumerate}
\item [(i)] $ R(x,h^0(x))=0$,
\item [(ii)] ${\rm det} \left(\displaystyle \frac{\partial
R_i}{\partial y_j}(x,h^{0}(x))\right)_{1\leq i,j\le r }\not
=0.$
\end{enumerate} Then there
exists a positive integer
$k=k(h^0)$ such that the following holds. If
$h:(\C^{q},0)\to (\C^r,0)$ is a formal map
such that $R(x,h(x))=0$ and $j_0^{k}h=j_0^{k}h^0$,
then necessarily $h(x)=h^{0}(x)$.
\end{Lem} \begin{proof} 
We may write
\begin{equation}\Label{eq56}
R(x,y)-R(x,t)=P(x,y,t)\cdot(y-t)
\end{equation} where $P$ is an $r\times r$
matrix with entries in $\C \dbl x,y,t\dbr$ satisfying
$P(x,y,y)=\displaystyle \frac{\partial R}{\partial
y}(x,y)$. By assumption, we know that 
${\rm det\,} P(x,h^0(x),h^0(x))\not = 0$.  This
implies that one can find an integer $k$ such
that if $h:(\C^q,0)\to (\C^r,0)$ is a formal mapping which
agrees up to order $k$ with
$h^0$, then
${\rm det\,}P(x,h^0(x),h(x))\not = 0$. If, in addition, $h$
satisfies $R(x,h(x))=0$, it follows from
(\ref{eq56}) that 
$P(x,h^0(x),h(x))\cdot(h^0(x)-h(x))= 0$ in
$\C[[x]]$.  Since
${\rm det\,}P(x,h^0(x),h(x))\not =0$,
we conclude that $h(x)=h^0(x)$ and hence the
lemma follows.\end{proof}

\begin{proof}[Proof of Proposition {\rm \ref{fdhol}}]
First observe that if $H,H^0:(\C^N,0)\to (\C^{N'},0)$ 
are two formal mappings with ${\mathcal I}^H={\mathcal
I}^{H^0}$, then by Proposition \ref{both} (i) and in
view of (\ref{game}), necessarily for any
$\alpha\in
\N^{n'}$,
$q_{\alpha}\circ H=q_{\alpha}\circ H^0$. Since $\scrM'$ is
holomorphically nondegenerate, we may choose
$\alpha^1,\ldots,\alpha^{N'}\in
\N^{N'}$ and
$j_1,\ldots,j_{N'}\in \{1,\ldots,d'\}$ as in
Lemma \ref{gorelose}. For any $l=1,\ldots, N'$, we define a
formal map $R_{l}:(\C^{N}\times \C^{N'},0)\to (\C,0)$ as
follows
\begin{equation} 
R_{l}(Z,Z'):=q_{j_l,\alpha^l}(Z')-q_{j_l,\alpha^l}(H^0(Z)).
\end{equation}
Observe that $R_l(Z,H^0(Z))=0$, for $l=1,\ldots,N'$, and
moreover, since $\Rk H^0=N'$, by (\ref{det}) and e.g.
Proposition 5.3.5 in \cite{BER}, we have
\begin{equation}\Label{linda}
{\rm det}\left(\frac{\partial q_{j_l,\alpha^{l}}}{\partial
Z'_m}(H^0(Z))\right)_{1\le l,m\le N'}\not
= 0,
\end{equation}
or equivalently, 
$${\rm det}\left(\frac{\partial R_l}{\partial
Z'_m}(Z,H^0(Z))\right)_{1\le l,m\le N'}\not
= 0.$$
 By
Lemma
\ref{prop22}, there exists a positive integer
$k=k(H^0)$ such that if $H:(\C^N,0)\to (\C^{N'},0)$ is a
formal map satisfying $R_l(Z,H(Z))=0$, for $l=1,\ldots,N'$,
and $j_0^kH=j_0^kH^{0}$, then $H=H^0$. On the other hand,
as mentioned in the beginning of the proof, if ${\mathcal
I}^H={\mathcal I}^{H^0}$, then
$R_{l}(Z,H(Z))=R_{l}(Z,H^0(Z))=0$, for $l=1,\ldots,N'$. This
completes the proof of Proposition \ref{fdhol}.
\end{proof}

\begin{proof}[Proof of Proposition {\rm \ref{link}}]
Since $M'$ is real-analytic, we may assume that the
corresponding formal mappings $\rho'(Z',\zeta')$ and $q_{\alpha}(Z')$
given in (\ref{game}) are convergent. First, note if there exists a
convergent  map $\check H:(\C^N,0)\to
(\C^{N'},0)$ such that ${\mathcal I}^{H}={\mathcal
I}^{\check H}$, then, since  ${\mathcal I}^{\check H}$ is
convergent, so is ${\mathcal
I}^{H}$. Now, assume that
${\mathcal I}^H$ is convergent. By Proposition
\ref{both} (ii) and in view of (\ref{game}), 
\begin{equation}\Label{artine}
r_{\alpha}(Z):=q_{\alpha}(H(Z))
\end{equation}
is a convergent mapping for all
$\alpha \in \N^{n'}$. By Artin's approximation theorem
\cite{artin1}, for any positive integer $\kappa$, there
exists a convergent map
$H^{\kappa}:(\C^N,0)\to (\C^{N'},0)$ which agrees with $H$
up to order $\kappa$ and such that
$q_{\alpha}(H^{\kappa}(Z))=r_{\alpha}(Z)$, for all $\alpha
\in
\N^{n'}$. It follows from (\ref{game}) and
(\ref{artine}) that 
$\rho'(H^{\kappa}(Z),\zeta')=\rho'(H(Z),\zeta')$ and hence
$${\mathcal
I}^{H}=(\rho'(H(Z),\zeta'))=(\rho'(H^{\kappa}(Z),\zeta'))={\mathcal
I}^{H^{\kappa}}.$$ This completes the proof of Proposition \ref{link}
in the convergent case. In the case when $M'$ is real-algebraic, the
$q_{\alpha}(Z')$ given by (\ref{game}) are algebraic. As
before, if there exists an algebraic map $\check H:(\C^N,0)\to
(\C^{N'},0)$ such that ${\mathcal I}^{H}={\mathcal
I}^{\check H}$, then, since  ${\mathcal I}^{\check H}$ is
algebraic, so is ${\mathcal
I}^{H}$. Moreover, it follows
from the algebraic version of Artin's theorem \cite{artin2}
that, in this case, one can choose
$H^{\kappa}$ as above to be algebraic so that
${\mathcal I}^{H}={\mathcal I}^{H^{\kappa}}$. The proof of
the proposition is now complete.
\end{proof}

For the proof of Proposition \ref{convhol}, we need the
following lemma whose proof is in the spirit of that of
Lemma \ref{prop22} but also makes use of Artin's
approximation theorem \cite{artin1}. We refer the reader to 
Proposition 4.2 of \cite{Mirh} for the proof of this lemma.

\begin{Lem}\Label{MRL} Let
$R(x,y)=(R_1(x,y),\ldots,R_r(x,y))\in (\C\{x,y\})^r$, $x\in
\C^q$, $y\in \C^r$, and 
$h:(\C^q,0)\to (\C^r,0)$  a formal map satisfying 
$R(x,h(x))=0$.  If ${\rm det}\big(\displaystyle
\frac{\partial R}{\partial y}(x,h(x))\big)\not
= 0$ in $\C[[x]]$, then $h(x)$ is convergent. 
\end{Lem}

\begin{proof}[Proof of Proposition {\rm \ref{convhol}}]
By Proposition \ref{both} (ii), if ${\mathcal I}^H$ is
convergent, then, in view of (\ref{game}), it follows that
for any $\alpha\in \N^{n'}$ and $j=1,\ldots,d'$, the formal
power series 
$r_{j,\alpha}(Z):=q_{j,\alpha}(H(Z))$ is convergent. Since
$M'$ is holomorphically nondegenerate, we may choose
$\alpha^1,\ldots,\alpha^{N'}\in
\N^{N'}$ and
$j_1,\ldots,j_{N'}\in \{1,\ldots,d'\}$ as in
Lemma \ref{gorelose}. For any $l=1,\ldots, N'$, we define a
convergent map $R_{l}:(\C^{N}\times \C^{N'},0)\to (\C,0)$ as
follows
\begin{equation} 
R_{l}(Z,Z'):=q_{j_l,\alpha^l}(Z')-r_{j_l,\alpha^l}(Z).
\end{equation}
Observe that $R_l(Z,H(Z))=0$, $l=1,\ldots,N'$, and
moreover, since $\Rk H=N'$, by (\ref{det}) and e.g.
Proposition 5.3.5 in \cite{BER}, we have 
\begin{equation}\Label{linda'}
{\rm det}\left(\frac{\partial q_{j_l,\alpha^{l}}}{\partial
Z'_m}(H(Z))\right)_{1\le l,m\le N'}\not
= 0,
\end{equation}
or equivalently, 
$${\rm det}\left(\frac{\partial R_l}{\partial
Z'_m}(Z,H(Z))\right)_{1\le l,m\le N'}\not
= 0.$$
We may now apply Lemma \ref{MRL} to conclude that $H$ is
convergent. The proof of Proposition \ref{convhol} is
complete.\end{proof}

\begin{proof}[Proof of Theorems {\rm \ref{plus1}} and {\rm
\ref{plus2}}] Theorem  \ref{plus1} is a
consequence of Theorem \ref{fdi} and Proposition
\ref{fdhol}, while Theorem  \ref{plus2} follows from
Theorem
\ref{convideal} and Proposition
\ref{convhol}. 
\end{proof}

\section{Proofs of Theorems \ref{fd},
\ref{conv},
\ref{approx1}, \ref{approx2}, \ref{plus3} and Corollaries
\ref{hypalg} and \ref{inj}}\Label{saved}

We begin with the following lemma, which will be used in
the proofs in this section.

\begin{Lem}\Label{finite}
Let $\scrM, \scrM'\subset \C^N\times \C^{N}$ be two formal
generic manifolds of the same codimension $d$ and
$H:(\C^N,0)\to (\C^N,0)$ a formal finite map. Then $\Rk
H=N$. Moreover, if the complexification
${\mathcal H}$ of $H$ maps $\scrM$ into $\scrM'$, then $H$
is not totally degenerate.
\end{Lem}
\begin{proof} The proof that $\Rk H=N$ is standard (see
e.g.\ Theorem 5.1.37 of \cite{BER}).  To prove the second
part of the lemma, it suffices to show that if
$\gamma(\zeta,t)$ is a Segre variety mapping as defined in
(\ref{gamma}) relative to $\scrM$, then $\Rk (H\circ
v^1)=n$, where 
$n=N-d$ and $v^1(t)=\gamma(0,t)$ as in (\ref{vj}). We
claim that the formal map
$ H\circ v^1$ is finite. Indeed, it
is a composition of the
finite map $H$ and of the formal map $v^1$ whose rank at 0
is
$n$ and hence is finite. The claim follows from the fact
that the composition of two formal finite mappings is again
finite. (This could be seen by e.g.\ making use of
Proposition 5.1.5 of \cite{BER}.) As before, the fact
that $H\circ v^1$ is finite implies that $\Rk (H\circ
v^1)=n$, which completes the proof of the lemma.
\end{proof}

\begin{proof}[Proof of Theorem {\rm \ref{fd}}] Without
loss of generality, we may assume that $p=p'=0$. Since
$M,M'\subset \C^N$ are smooth generic submanifolds through
the origin,  we can consider the associated formal
generic manifolds $\scrM,\scrM'\subset \C^N\times \C^N$
as described in \S \ref{ideal}. In this case, the
complexification ${\mathcal H}$ of any formal map
$H:(\C^N,0)\to (\C^N,0)$ sending $M$ into $M'$ sends
$\scrM$ into $\scrM'$. Since the given formal map
$H^0$ is finite, it follows from Lemma
\ref{finite} that $H^0$ is not totally degenerate and
$\Rk H^0=N$. Theorem \ref{fd} is then a consequence of
Theorem \ref{plus1}.
\end{proof}

\begin{proof}[Proof of Theorem {\rm \ref{conv}}]
Without loss of generality, we may assume that $p=p'=0$.
Since the given formal map $H$ is finite, it follows from Lemma
\ref{finite} that $H$ is not totally degenerate and
$\Rk H=N$. Theorem \ref{conv} is then a
consequence of Theorem \ref{plus2}.\end{proof}

For the proofs of Theorems \ref{approx1} and \ref{approx2},
we need the following lemma.

\begin{Lem}\Label{mir}
Let $I\subset \C \dbl Z,\zeta\dbr$ and $J\subset \C\dbl
Z',\zeta'\dbr$ be two ideals and $H,\check H:(\C^N_Z,0)\to
(\C^{N'}_{Z'},0)$ be two formal mappings. Let ${\mathcal
H},{\check \mathcal H}$ be the complexifications of $H$ and
$\check H$ respectively as defined in {\rm (\ref{H})}.
Assume that:
\begin{enumerate}
\item [(i)] $J$ is a real ideal;
\item [(ii)] $J\subset {\mathcal H}_{*}(I)$, where
${\mathcal H}_{*}(I)$, the pushforward of $I$ by
${\mathcal H}$ as defined by {\rm (\ref{push})};
\item [(iii)] $J^{\check H}\subset J^{H}$, where the ideals
$J^{\check H},J^{ H}\subset \C \dbl Z,\zeta'\dbr$ are
defined by {\rm (\ref{ballot})}.
\end{enumerate}
Then $J\subset \check {\mathcal H}_{*}(I)$.
\end{Lem}
\begin{proof} Let $s_1(Z',\zeta'),\ldots, s_{m}(Z',\zeta')$
be generators of $J$ in $\C \dbl Z',\zeta'\dbr$. As usual,
we write $s(Z',\zeta')=(s_1(Z',\zeta'),\ldots,
s_{m}(Z',\zeta'))$ and
$J=(s(Z',\zeta'))$. We set
\begin{equation}\Label{small}
\tilde
s(Z',\zeta'):=\overline{s}(\zeta',Z').
\end{equation}
 By the reality of $J$,
it follows that we also have $J=(\tilde s(Z',\zeta'))$.
Hence there exists an $m\times m$ matrix with entries in 
$\C \dbl Z',\zeta'\dbr$ such that 
\begin{equation}\Label{big}
s(Z',\zeta')=u(Z',\zeta') \tilde s(Z',\zeta').
\end{equation}
Note that in view of  (\ref{ballot}), the ideals
$J^{H}$ and $J^{\check H}$ are generated by the components
of
$s(H(Z),\zeta')$ and $s(\check H(Z),\zeta')$ respectively
in $\C\dbl Z,\zeta'\dbr$. Hence, by the inclusion (iii), we
have
\begin{equation}\Label{doubt}
s(\check H(Z),\zeta')=a(Z,\zeta')s(H(Z),\zeta'),
\end{equation}
where $a(Z,\zeta')$ is an $m\times m$ matrix with entries in
$\C\dbl Z,\zeta'\dbr$. By taking complex conjugates, it
follows from (\ref{doubt}) that we also have
\begin{equation}\Label{doubt2}
\overline{s}(\overline{\check
H}(\zeta),Z')=\bar{a}(\zeta,Z')\overline{s}(\bar{H}(\zeta),Z').
\end{equation}
To prove the lemma, we must show that the
components of $s(\check H(Z),\overline{\check H}(\zeta))$
are in $I$. For this, using (\ref{doubt}), (\ref{big}),
(\ref{small}) and (\ref{doubt2}), we have
\begin{equation}\Label{influence}
\begin{aligned}
s(\check H(Z),\overline{\check H}(\zeta))&=
a(Z,\overline{\check H}(\zeta))\, s(H(Z),\overline{\check
H}(\zeta))\\
&= a(Z,\overline{\check H}(\zeta))\, u(
H(Z),\overline{\check H}(\zeta))\,  \tilde s(
H(Z),\overline{\check H}(\zeta))\\
&=a(Z,\overline{\check H}(\zeta))\, u(
H(Z),\overline{\check H}(\zeta))\,
\overline{s}(\overline{\check H}(\zeta),H(Z))\\
&= a(Z,\overline{\check H}(\zeta))\, u(
H(Z),\overline{\check H}(\zeta))\,
\overline{a}(\zeta,H(Z))\overline{s}(\overline{H}(\zeta),H(Z))\\
&=a(Z,\overline{\check H}(\zeta))\, u(
H(Z),\overline{\check H}(\zeta))\,
\overline{a}(\zeta,H(Z))\tilde{s}(H(Z),\overline{H}(\zeta)).
\end{aligned}
\end{equation}
By (ii), the components of
$\tilde{s}(H(Z),\overline{H}(\zeta))$ are in $I$ and hence,
by (\ref{influence}), so are the components of $s(\check
H(Z),\overline{\check H}(\zeta))$. The proof of the lemma
is complete.
\end{proof}

The following lemma, which is a reformulation of an
observation due to Merker, is an immediate consequence of
Lemma 
\ref{mir}.

\begin{Lem}\Label{merker}
Let $\scrM\subset \C^N\times \C^N$ and $\scrM'\subset
\C^{N'}\times \C^{N'}$ be two formal generic manifolds and
$H,\check H:(\C^{N},0)\to (\C^{N'},0)$ be two formal
mappings whose complexifications are denoted by ${\mathcal
H}$ and  $\check {\mathcal H}$ respectively. Assume that
${\mathcal H}$ sends $\scrM$ into $\scrM'$ and that the
reflection ideals ${\mathcal I}^H$ and ${\mathcal
I}^{\check H}$ are the same. Then 
$\check {\mathcal H}$ also sends
$\scrM$ into $\scrM'$. 
\end{Lem}

\begin{proof}[Proof of Theorem {\rm \ref{plus3}}]
Since $M$ is of finite type at 0 and the formal map $H$ is not totally
degenerate, by Theorem
\ref{convideal}, the reflection ideal ${\mathcal I}^H$ is
convergent. By Proposition \ref{link}, for any
positive integer $\kappa$, there exists a convergent map
$H^{\kappa}:(\C^N,0)\to (\C^N,0)$ which agrees with $H$ up
to order $\kappa$ such that
${\mathcal I}^{H}={\mathcal I}^{H^{\kappa}}$. By Lemma
\ref{merker}, it follows that ${\mathcal
H}^{\kappa}$ maps $\scrM$ into $\scrM'$ and
hence $H^{\kappa}$ maps
$M$ into
$M'$. The proof of Theorem \ref{plus3} is complete.
\end{proof}

\begin{proof}[Proof of Theorem {\rm \ref{approx1}}] 
Without loss of generality, we may assume that $p=p'=0$. 
Since the
given formal map $H$ is finite, it follows from Lemma
\ref{finite} that $H$ is not totally degenerate.
Theorem \ref{approx1} is then a consequence of Theorem
\ref{plus3}.
\end{proof}

\begin{proof}[Proof of Theorem {\rm \ref{approx2}}] 
Without loss of generality, we may assume that $p=p'=0$. 
Since $M$ is connected and of finite type at some point, by
Lemma 13.3.2 of \cite{BER}, there is no germ of a
nonconstant holomorphic function
$h:(\C^N,0)\to \C$ with $h(M)\subset \bR$. It follows from 
Theorem \ref{algideal} that the reflection ideal ${\mathcal
I}^H$ of the given local holomorphic map $H$ is algebraic. By
Proposition
\ref{link}, for any positive integer $\kappa$, there exists an
algebraic map
$H^{\kappa}:(\C^N,0)\to (\C^N,0)$ which agrees with $H$ up
to order $\kappa$ such that
${\mathcal I}^{H}={\mathcal I}^{H^{\kappa}}$. By Lemma
\ref{merker}, it follows that ${\mathcal
H}^{\kappa}$ maps $\scrM$ into $\scrM'$ and
hence $H^{\kappa}$ maps
$M$ into
$M'$. The proof of Theorem \ref{approx2} is complete.
\end{proof}

\begin{proof}[Proof of Corollary {\rm \ref{hypalg}}]
Let $(M,p)$ and $(M',p')$ be two germs of
biholomorphically equivalent real-algebraic hypersurfaces in
$\C^N$. If there is no point  of finite type in $M$
arbitrarily close to $p$, then $(M,p)$ is Levi-flat and so
is $(M',p')$. Hence both $(M,p)$ and $(M',p')$ are
algebraically equivalent to a real hyperplane in $\C^N$. If
$M$ contains points of finite type arbitrarily close to $p$,
then we may apply Corollary \ref{alg} to conclude that
$(M,p)$ and $(M',p')$ are algebraically equivalent. The
proof of Corollary \ref{hypalg} is complete.
\end{proof}

\begin{proof}[Proof of Corollary {\rm \ref{inj}}]
By Theorem \ref{fd} with $(M,p)=(M',p')$ and $H^0={\rm Id}$,
the identity map of $(\C^N,p)$, there exists a positive
integer
$K$ such that if $H:(\C^N,p)\to (\C^N,p)$ is a formal map
sending $M$ into itself with $j_p^KH=j_p^K{\rm Id}$, then
 $H={\rm Id}$. Let $H^1,H^2:(\C^N,p)\to (\C^N,p)$ be two
invertible formal mappings sending $M$ into itself and such
that $j_p^KH^1=j_p^KH^2$. If $H=H^1\circ (H^2)^{-1}$, then
$H$ is formal map sending $(M,p)$ into itself such that
$j_p^KH=j_p^K{\rm Id}$. Hence $H=H^1\circ (H^2)^{-1}={\rm
Id}$, i.e.\ $H^1=H^2$. The second part of Corollary
\ref{inj} is an immediate application of Theorem \ref{conv}.
\end{proof}

\section{Remarks and open problems} \Label{open}
As mentioned in \S \ref{int},
holomorphic nondegeneracy is necessary for the conclusions of
Theorems \ref{fd} and \ref{conv} to hold. Indeed, if $(M,p)$ is a
germ of a smooth generic submanifold in
$\C^N$ which is holomorphically degenerate at $p$, then for any
positive integer $K$, there exist a formal invertible mapping
$H:(\C^N,p)\to (\C^N,p)$ sending $M$ into itself and agreeing with
the identity map Id up to order $K$ at $p$ but such that $H\not =
{\rm Id}$ (see \cite{CAG} Theorem 3 and \cite{MA} Theorem 2.2.1).
Similarly, if $(M,p)$ is a
germ of a real-analytic generic submanifold in
$\C^N$ which is holomorphically degenerate at $p$, then there exist
(infinitely many) nonconvergent formal invertible self-mappings of
$(M,p)$ (see
\cite{Asian}).

In constrast to holomorphic nondegeneracy, the finite type condition
in Theorems \ref{fd} and \ref{conv} does not seem to be necessary.
More precisely, we conjecture the following. If $M\subset \C^N$ is a
connected
 holomorphically nondegenerate real-analytic generic submanifold of
finite type at some point, then for any $p\in M$, ${\rm
Aut}\,(M,p)={\mathcal F}(M,p)$. Here, we recall that ${\rm
Aut}\,(M,p)$ is the stability group of $(M,p)$ and ${\mathcal F}(M,p)$
is the group of formal invertible self-mappings of $(M,p)$. We also
conjecture that if $M$ is as above, then
for every $p\in M$, there exists a positive integer $K=K(p)$ such
that the jet mapping $j_p^K:{\rm Aut}\, (M,p)\to G^K(\C^N,p)$ is
injective, where $G^K(\C^N,p)$ is the jet group of order $K$ at $p$.
It follows from Corollary \ref{inj} that the above
conjectures hold for all points $p$ in a Zariski open subset of $M$.

Another question concerning the structure of ${\mathcal F}(M,p)$ is
the following. Under the assumptions of Corollary
\ref{inj}, is the image of the group homomorphism
$j_p^K: \mathcal F(M,p) \to G^K(\C^N,p)$ a
closed Lie subgroup of the jet group $G^K(\C^N,p)$, for some suitable
integer $K$? The question is open even when $M$ is
real-analytic, in which case $\mathcal F(M,p)={\rm Aut}\, (M,p)$, by
Corollary \ref{inj}. It is known that the answer is positive 
if
$M$ is finitely nondegenerate  and of finite type at $p$ (see
\cite{Asian} for the hypersurface case and \cite {Z1} for higher
codimension).  In fact it is shown in \cite{MA} 
that in this case the image is actually a totally real algebraic Lie
subgroup of $G^K(\C^N,p)$ for a precise value of $K$.

Finally, concerning algebraic equivalence, in view of Corollary
\ref{hypalg}, one is led to conjecture that biholomorphic equivalence
implies algebraic equivalence for germs of real-algebraic
submanifolds in $\C^N$. To the knowledge of the authors, the question
is still open even for germs of generic real-algebraic submanifolds
of codimension higher than one.

\end{document}